\documentclass[final,leqno]{siamltex}

\usepackage{times}
\usepackage{graphics,epsfig,graphicx,color}
\usepackage{amsmath}
\usepackage{mathrsfs,amsfonts}
\usepackage{dsfont}
\usepackage{url}

\newtheorem{prop}[theorem]{Proposition}

\newcommand{\bb}{{\bf b}}
\newcommand{\bd}{{\bf d}}
\newcommand{\be}{{\bf e}}
\newcommand{\bl}{{\bf l}}
\newcommand{\br}{{\bf r}}
\newcommand{\bs}{{\bf s}}
\newcommand{\brho}{{\boldsymbol{\rho}}}
\newcommand{\bu}{{\bf u}}

\newcommand{\bx}{{\bf x}}
\newcommand{\by}{{\bf y}}
\newcommand{\bz}{{\bf z}}

\newcommand{\bfdelta}{\boldsymbol{\delta}}

\newcommand{\bftheta}{\boldsymbol{\theta}}
\newcommand{\bflambda}{\boldsymbol{\lambda}}

\newcommand{\bfeta}{\boldsymbol{\eta}}
\newcommand{\bfrho}{\boldsymbol{\rho}}

\newcommand{\cB}{{\mathcal B}}
\newcommand{\cD}{{\mathcal D}}
\newcommand{\cE}{{\mathcal E}}
\newcommand{\cH}{{\mathcal H}}

\newcommand{\cJ}{{\mathcal J}}
\newcommand{\cK}{{\mathcal K}}
\newcommand{\cM}{{\mathcal M}}

\newcommand{\cQ}{{\mathcal Q}}
\newcommand{\cR}{{\mathcal R}}
\newcommand{\cF}{{\mathcal F}}
\newcommand{\CC}{{\mathbb{C}}}

\title{A model reduction approach to numerical inversion for a
  parabolic partial differential equation}

\author{
Liliana Borcea\footnotemark[1]
\and
Vladimir Druskin\footnotemark[2]
\and
Alexander V. Mamonov\footnotemark[3]
\and
Mikhail Zaslavsky\footnotemark[4]
}

\begin{document}

\maketitle

\renewcommand{\thefootnote}{\fnsymbol{footnote}}

\footnotetext[1]{Department of Mathematics, University of Michigan,
  2074 East Hall, 530 Church Street, Ann Arbor, MI 48109-1043
  (borcea@umich.edu)} \footnotetext[2]{Schlumberger-Doll Research
  Center, 1 Hampshire St., Cambridge, MA 02139-1578
  (druskin1@slb.com)} \footnotetext[3]{Schlumberger, 
  3750 Briarpark Dr., Houston, TX 77042
  (amamonov@slb.com)} \footnotetext[4]{Schlumberger-Doll Research
  Center, 1 Hampshire St., Cambridge, MA 02139-1578
  (mzaslavsky@slb.com)}

\begin{abstract}
We propose a novel numerical inversion algorithm for the coefficients 
of parabolic partial differential equations, based on model reduction. 
The study is
motivated by the application of controlled source electromagnetic
exploration, where the unknown is the subsurface electrical
resistivity and the data are time resolved surface measurements of the
magnetic field. The algorithm presented in this paper considers
inversion in one and two dimensions. The reduced model is obtained
with rational interpolation in the frequency (Laplace) domain and a
rational Krylov subspace projection method. It amounts to a non-linear
mapping from the function space of the unknown resistivity to the
small dimensional space of the parameters of the reduced model. We use
this mapping as a non-linear preconditioner for the Gauss-Newton
iterative solution of the inverse problem. The advantage of the
inversion algorithm is twofold.  First, the non-linear preconditioner
resolves most of the nonlinearity of the problem.  Thus the iterations
are less likely to get stuck in local minima and the convergence is
fast. Second, the inversion is computationally efficient because it
avoids repeated accurate simulations of the time-domain response.  
We study the stability of the inversion algorithm
for various rational Krylov subspaces, and assess its performance with
numerical experiments.
\end{abstract}



\pagestyle{myheadings}
\thispagestyle{plain}

\section{Introduction}
\label{sec:intro}
Inverse problems for parabolic partial differential equations arise in
applications such as groundwater flow, solute transport and controlled
source electromagnetic oil and gas exploration.  We consider the
latter problem, where the unknown $r({\bf x})$ is the electrical
resistivity, the coefficient in the diffusion Maxwell system satisfied
by the magnetic field ${\bf H}(t,{\bf x})$
\begin{equation}
\label{eq:2}
-\nabla \times [ r({\bf x}) \nabla \times {\bf H}(t,{\bf x})] = 
\frac{\partial {\bf H}(t,{\bf x})}{\partial t}, 
\end{equation} 
for time $t > 0$ and ${\bf x}$ in some spatial domain. The data 
from which $r(\bx)$ is to be determined are the
time resolved measurements of ${\bf H}(t,{\bf x})$ at receivers
located on the boundary of the domain.

Determining $r(\bx)$ from the boundary measurements is challenging 
especially because the problem is ill-posed and thus sensitive to noise. 
A typical numerical approach is to minimize a functional given by 
the least squares data misfit and a regularization term, using 
Gauss-Newton or non-linear conjugate gradient methods 
\cite{newman2000three, zhdanov, oristaglio1994inversion}.  
There are two main drawbacks. First, the
functional to be minimized is not convex and the optimization
algorithms can get stuck in local minima.  The lack of convexity can
be overcome to some extent by adding more regularization at the cost
of artifacts in the solution. Nevertheless, convergence may be very
slow \cite{oristaglio1994inversion}.  Second, evaluations of the
functional and its derivatives are computationally expensive, because
they involve multiple numerical solutions of the forward problem. In
applications the computational domains may be large with meshes
refined near sources, receivers and regions of strong
heterogeneity. This results in a large number of unknowns in the
forward problem, and time stepping with such large systems is
expensive over long time intervals.


\begin{figure}[ht!]
\begin{center}
\includegraphics[width=0.4\textwidth]{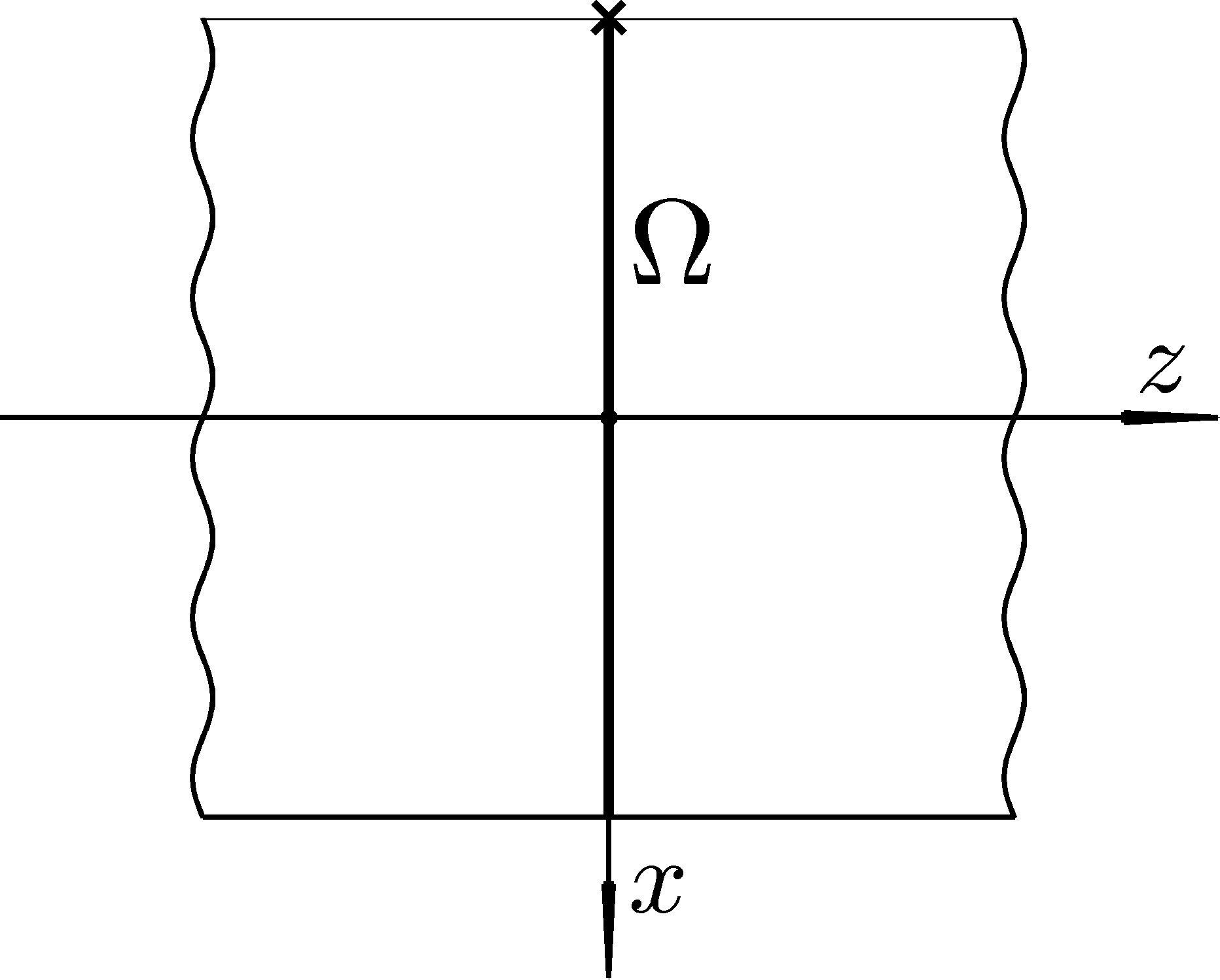} \hskip0.1\textwidth
\includegraphics[width=0.4\textwidth]{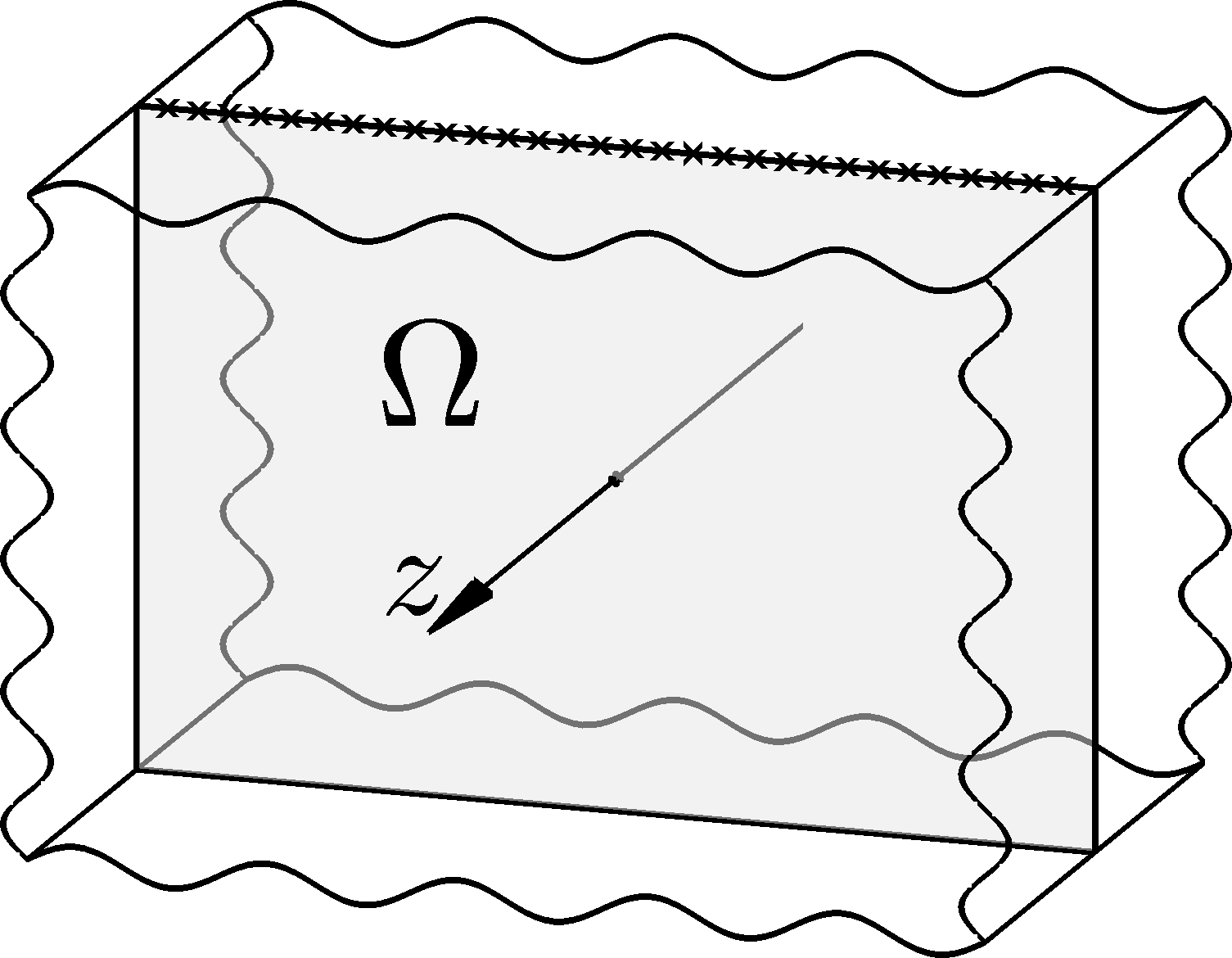}
\caption{Examples of spatial domains for equation (\ref{eq:2}) in 
$\mathbb{R}^2$ (left) and $\mathbb{R}^3$ (right). The medium is constant
in the direction $z$, which is transversal to $\Omega$. The accessible
boundary $\cB_A \subset \cB = \partial \Omega$ is marked with $\times$.
}
\label{fig:setting}
\end{center}
\end{figure}

We propose a numerical inversion approach that arises when considering
the inverse problem in the model reduction framework. We consider one
and two dimensional media, and denote the spatial variable by
${\bf x} = (x,z)$, with $z \in \mathbb{R}$ and $x \in
\Omega$, a simply connected domain in $\mathbb{R}^{n}$ with boundary
$\mathcal{B}$, for $n = 1,2$. The setting is illustrated in figure
\ref{fig:setting}.

The resistivity is $r = r(x)$, and assuming that
\begin{equation}
{\bf H} = u(t,x) {\bf e}_z, 
\end{equation}
we obtain from (\ref{eq:2}) the parabolic equation 
\begin{equation}
\nabla \cdot \left[ r(x) \nabla u(t,x) \right] = 
\frac{\partial u(t, x)}{\partial t}, 
\label{eq:3d}
\end{equation}
for $t > 0$ and $x \in \Omega$. The boundary $\mathcal{B}$ 
consists of an accessible part $\mathcal{B}_A$, which supports
the receivers at which we make the measurements and the initial excitation
\begin{equation}
u(0,x) = u_o(x), \qquad {\rm supp}\{u_o\} \subset \cB_A,
\label{eq:5d} 
\end{equation} 
and an inaccessible part $\cB_I = \cB \setminus \cB_A$ where we set
\begin{equation}
\label{eq:4dd}
 u(t, x) = 0, \quad {x} \in \cB_I.
\end{equation}
The boundary condition at $\mathcal{B}_A$ is 
\begin{equation}
\label{eq:4d}
{\bf n}(x) \cdot \nabla u(t,x) = 0, \quad x \in \mathcal{B}_A,
\end{equation}
where ${\bf n}$ is the outer normal.  We choose the boundary
conditions (\ref{eq:4dd}-\ref{eq:4d}) to simplify the
presentation. Other (possibly inhomogeneous) conditions can
be taken into account and amount to minor modifications of the reduced
models described in this paper.

The sources and receivers, are located on 
the accessible boundary $\mathcal{B}_A$, and provide 
knowledge of the operator
\begin{equation}
\cM(u_o) = \left. u(t,x) \right|_{x \in \cB_A}, \quad t > 0,
\label{eqn:contmeas}
\end{equation}
for every $u_o$ such that ${\rm supp} \{ u_o \} \subset \cB_A$.
Here $u(t, x)$ is the solution of (\ref{eq:3d}) with the initial
condition (\ref{eq:5d}). The inverse problem is to determine
the resistivity $r(x)$ for $x \in \Omega$ from $\cM$. 
Note that the operator $\cM$ is highly non-linear in $r$, but it is linear 
in $u_o$. This implies that in one dimension where  $\Omega$ is an 
interval, say $\Omega=(0,1)$, and the accessible boundary is a 
single point $\mathcal{B}_A = \{x = 0\}$,  $\cM$
is completely defined by  $u_o(x) = \delta(x)$.  
All the information about $r(x)$ is contained in a single 
function of time
\begin{equation}
y(t) = \cM(\delta(x)) = u(t, 0), \quad t>0.
\label{eqn:y1Dcont}
\end{equation}

To ease the presentation, we begin by describing in detail the model
reduction inversion method for the one dimensional case. 
Then we show how
to extend it to the two dimensional case. The reduced model is
obtained with a rational Krylov subspace projection method. Such
methods have been applied to forward problems for parabolic equations
in \cite{borner2008fast,druskin2009solution,
frommer2008matrix}, and to inversion in one dimension in 
\cite{dsz2012} and multiple dimensions in \cite{beattietal}.
The reduced models in \cite{dsz2012} are $\mathcal{H}_2$ rational
interpolants of the transfer function defined by the Laplace transform
of $y(t)$ from (\ref{eqn:y1Dcont}). In this paper we build on the
 results in \cite{dsz2012} to study in more detail and improve the 
inversion method in one dimension, and to extend it to two dimensions.

The reduced order models allow us to replace the solution of the full
scale parabolic problem by its low order projection, thus resolving
the high computational cost inversion challenge mentioned above.
Conventionally, the reduced order model (the rational approximant of
the transfer function) is parametrized in terms of its poles and
residues. We parametrize it instead in terms of the coefficients of
its continued fraction expansion.  The rational approximant can be
viewed as the transfer function of a finite difference discretization
of (\ref{eq:3d}) on a special grid with rather few points, known in
the literature as the optimal or spectrally matched grid, or as a
finite-difference Gaussian rule \cite{druskin2000gaussian}.  The
continued fraction coefficients are the entries in the finite
difference operator, which are related to the discrete resistivities.

To mitigate the other inversion challenge mentioned above, we
introduce a non-linear mapping ${\mathcal \cR}$ from the function space
of the unknown resistivity $r$ to the low-dimensional space of the
discrete resistivities.  This map appears to resolve the nonlinearity
of the problem and we use it in a preconditioned Gauss-Newton
iteration that is less likely to get stuck in local minima and
converges quickly, even when the initial guess is far from the true
resistivity.

To precondition the problem, we map the measured data to the
discrete resistivities. The inverse problem is ill-posed, so we
limit the number of discrete resistivities (i.e. the size of the
reduced order model) computed from the data. This number depends on
the noise level and it is typically much less than the dimension of
models used in conventional algorithms, where it is determined by the
accuracy of the forward problem solution. This represents another
significant advantage of our approach.

The paper is organized as follows: We begin in section
\ref{sec:rominv} with a detailed description of our method in one
dimension. The inversion in two dimensional media is described in section
\ref{sec:2D}. Numerical results in one and two dimensions are
presented in section \ref{sec:numres}. We conclude with a summary in
section \ref{sec:conclusion}. The technical details of the computation
of the non-linear preconditioner $\cR$ and its Jacobian are given
in Appendix \ref{app:compjac}.

\section{Model order reduction for inversion in one dimensional media}
\label{sec:rominv}

In this section we study in detail the use of reduced order models for
inversion in one dimension. Many of the ideas presented here are used
in section \ref{sec:2D} for the two-dimensional inversion.  We begin
in section \ref{subsec:semidisc} by introducing a semi-discrete
analogue of the continuum inverse problem, where the differential
operator in $x$ in equation (\ref{eq:3d}) is replaced by a
matrix. This is done in any numerical inversion, and we use it from
the start to adhere to the conventional setting of model order
reduction, which is rooted in linear algebra. The projection-based
reduced order models are described in section
\ref{subsec:projrom}. They can be parametrized in terms of the
coefficients of a certain reduced order finite difference scheme, used
in section \ref{subsec:finitediff} to define a pair of non-linear
mappings.  They define the objective function for the non-linearly
preconditioned optimization problem, as explained in section
\ref{subsec:optim}. To give some intuition to the preconditioning
effect, we relate the mappings to so-called optimal grids in section
\ref{subsec:optgrid}. The stability of the inversion is addressed in
section \ref{subsec:ratinterp}, and our regularization scheme is given
in section \ref{subsec:regular}. The detailed description of the one
dimensional inversion algorithm is in section \ref{subsec:invalg}.

\subsection{Semi-discrete inverse problem}
\label{subsec:semidisc}

In one dimension the domain is an interval, which we scale to $\Omega
= (0,1)$, with accessible and inaccessible boundaries $\mathcal{B}_A =
\{0\}$ and $\mathcal{B}_I = \{1\}$ respectively. To adhere to the
conventional setting of model order reduction, we consider the
semi-discretized equation (\ref{eq:3d})
\begin{equation}
\frac{\partial {\bf u}(t)}{\partial t} = A(\br) {\bf u}(t),
\label{eqn:semidfwd}
\end{equation}
where 
\begin{equation}
A(\br) = - D^T \mbox{diag}(\br) D,
\label{eqn:Ar}                                                   
\end{equation}
is a symmetric and negative definite matrix, the discretization of
$\partial_x [ r(x) \partial _x]$.  The vector $\br \in \mathbb{R}^N_+$
contains the discrete values of $r(x)$ and the matrix $D \in
\mathbb{R}^{N \times N}$ arises in the finite difference
discretization of the derivative in $x$, for boundary conditions
(\ref{eq:4d}-\ref{eq:4dd}). The discretization is on a very fine
uniform grid with $N$ points in the interval $[0,1]$, and spacing $h =
1/(N+1)$. Note that our results do not depend on the dimension $N$,
and all the derivations can be carried out for a continuum
differential operator.  We let $A$ be a matrix to avoid unnecessary
technicalities.  The vector ${\bf u}(t) \in \mathbb{R}^N$ is the
discretization of $u(t,x)$, and the initial condition is
\begin{equation}
\bu(0) = \frac{1}{h} {\bf e}_1,
\label{eqn:semidfwdini}
\end{equation}
the approximation of a point source excitation at $\mathcal{B}_A$. The
time-domain response of the semi-discrete dynamical system
(\ref{eqn:semidfwd}) is given by
\begin{equation}
y(t; \br) = {\bf e}_1^T \bu(t),
\label{eqn:semidresp}
\end{equation}
where $\be_1 = (1, 0, 0, \ldots, 0)^T$. This is a direct analogue of
(\ref{eqn:y1Dcont}), i.e. it corresponds to measuring the solution of
(\ref{eqn:semidfwd}) at the left end point $\mathcal{B}_A = \{0\}$ of
the interval $\Omega$. We emphasize in the notation that the response
depends on the vector $\br$ of discrete resistivities.
 
We use (\ref{eqn:semidresp}) to define the \emph{forward map}
\begin{equation}
\cF : \mathbb{R}^N_+ \to C(0, +\infty)
\end{equation}
that takes the vector $\br \in \mathbb{R}^N_+$ to the time domain
response
\begin{equation}
\cF(\br) = y(\,\cdot\,; \br).
\label{eqn:forwardmap}
\end{equation}
The measured time domain data is denoted by
\begin{equation}
d(t) = y(t; \br^{\mbox{\tiny true}}) + \mathcal{N}(t), \quad t>0,
\label{eqn:semiddata}
\end{equation}
where $\br^{\mbox{\tiny true}}$ is the true (unknown) resistivity
vector and $\mathcal{N}(t)$ is the contribution of the measurements
noise and the discretization errors. The inverse problem is: {Given
  data $d(t)$ for $t \in [0, \infty)$, recover the resistivity
    vector.}


\subsection{Projection-based model order reduction}
\label{subsec:projrom}

In order to apply the theory of model order reduction we treat 
(\ref{eqn:semidfwd}-\ref{eqn:semidresp}) as a dynamical system with 
the time domain response 
\begin{equation}
y(t;\br) = {\bf e}_1^T e^{A(\br)t}\, \frac{{\bf e}_1}{h} = \bb^T
e^{A(r)t} \, \bb
\label{eqn:symy}
\end{equation} 
written in a symmetrized form using the source/measurement vector
\begin{equation}
\bb = {\bf e}_1/{\sqrt{h}}.
\end{equation}

The transfer function of the dynamical system is the Laplace transform
of (\ref{eqn:symy})
\begin{equation}
Y(s; \br) = \int_{0}^{+\infty} y(t; \br) e^{-st}dt = \bb^T (sI-A(\br))^{-1}
\bb, \quad s>0.
\label{eqn:Gtransfer}
\end{equation}
Since $A(\br)$ is negative definite, all the poles of the
transfer function are negative and the dynamical system is stable.

In model order reduction we obtain a reduced model $(A_m, \bb_m)$ so that
its transfer function
\begin{equation}
\label{eqn:Gmtransfer}
Y_m(s) = \bb_m^T (sI_m - A_m)^{-1} \bb_m
\end{equation}
is a good approximation of $Y(s;\br)$ as a function of $s$ in some
norm. Here $I_m$ is the $m \times m$ identity matrix, $A_m \in
\mathbb{R}^{m \times m}$, $\bb_m \in \mathbb{R}^m$ and $m \ll N$.
Note that while the matrix $A(\br)$ given by (\ref{eqn:Ar}) is sparse,
the reduced order matrix $A_m$ is typically dense. Thus, it has no
straightforward interpretation as a discretization of the differential
operator $\partial_x [ r(x) \partial_x]$ on some coarse grid with $m
$ points.

Projection-based methods search for reduced models of the form
\begin{equation} 
A_m = V^T A V, \quad \bb_m = V^T \bb, \quad V^T V = I_m,
\label{eqn:rom-project}
\end{equation}
where the columns of $V \in \mathbb{R}^{N \times m}$ form an
orthonormal basis of an $m$-dimensional subspace of $\mathbb{R}^N$ on
which the system is projected.
The choice of $V$ depends on the {matching conditions} for $Y_m$
and $Y$. They prescribe the sense in which $Y_m$ approximates
$Y$. Here we consider moment matching at  {interpolation
  nodes} $s_j \in [0, +\infty)$ that may be distinct or coinciding,
\begin{equation}
\left. \frac{\partial^k Y_m}{\partial s^k} \right|_{s = s_j} =
\left. \frac{\partial^k Y}{\partial s^k} \right|_{s = s_j}, 
\quad k=0,\ldots,2M_j-1, \quad j=1,\ldots,l.
\label{eqn:matching}
\end{equation}
The multiplicity of a node $s_j$ is denoted by $M_j$, so at
non-coinciding nodes $M_j=1$. The reduced order transfer function $Y_m$ 
matches $Y$ and its derivatives up to the order $2M_j-1$, and the size 
$m$ of the reduced model is given by $m = \sum_{j=1}^l M_j$.

Note from (\ref{eqn:Gmtransfer}) that $Y_m(s)$ is a rational function
of $s$, with partial fraction representation
\begin{equation}
Y_m(s) = \sum_{j=1}^{m} \frac{c_j}{s + \theta_j}, \quad c_j > 0, \quad
\theta_j > 0.
\label{eqn:partfrac}
\end{equation}
Its poles $-\theta_j$ are the eigenvalues of $A_m$ and the residues
$c_j$ are defined in terms of the normalized eigenvectors ${\bf z}_j$,
\begin{equation}
c_j = (\bb_m^T {\bf z}_j)^2, \quad A_m {\bf z}_j + \theta_j {\bf z}_j
= 0, \quad \|{\bf z}_j\| = 1, \quad j=1,\ldots,m.
\label{eqn:thetac}
\end{equation}
Thus, (\ref{eqn:matching}) is a rational interpolation problem. It is 
known  \cite{grimme1997krylov} to be equivalent to 
the projection (\ref{eqn:rom-project}) when the columns of $V$ form an
orthogonal basis of the rational Krylov subspace
\begin{equation}
\mathcal K_m(\bs) = \mbox{span} \left\{ (s_j I - A)^{-k} \bb
\; | \; j=1,\ldots,l; \; k=1,\ldots,M_j \right\}.
\label{eqn:ratkrylov}
\end{equation}
The interpolation nodes, obviously, should be chosen in the resolvent
set of $A(\br)$. Moreover, since in reality we need to solve a
limiting continuum problem, the nodes should be in the closure of the
intersection of the resolvent sets of any sequence of 
finite-difference operators that converge to the continuum
problem. This set includes $\CC\setminus (-\infty,0)$ for problems on
bounded domains. In our computations the interpolation nodes lie
on the positive real axis, since they correspond to working with the
Laplace transform of the time domain response.

Ideally, in the solution of the inverse problem we would like to
minimize the time (or frequency) domain data misfit in a quadratic
norm weighted in accordance with the statistical distribution of the
measurement error. When considering the reduced order model, it is
natural to choose the interpolation nodes that give the most accurate
approximation in that norm. Such interpolation is known in control
theory as $\cH_2$ (Hardy space) optimal, and in many cases the optimal
interpolation nodes can be found numerically
\cite{gugercin2008mathcal}. Moreover, it was shown in \cite{dsz2012},
that the solution of the inverse problem using reduced order models
with such interpolation nodes also minimizes the misfit functional (in
the absence of measurement errors). When such nodes are not available,
we can select some reasonable interpolation nodes chosen based on a
priori error estimates, for example, the so-called Zolotarev points or
their approximations obtained with the help of potential theory
\cite{druskin2009solution,zasletalgeoinv}. In most cases such choices
lead to interpolation nodes that are distributed geometrically.

\subsection{Finite difference parametrization of reduced order models}
\label{subsec:finitediff}

As we mentioned above, even though the reduced order model $(A_m,
\bb_m)$ comes from a finite difference operator $A$, it does not
retain its structure. In particular, $A_m$ is a dense matrix. The
model can be uniquely parametrized by $2m$ numbers, for example the
poles $-\theta_j$ and residues $c_j$ in the representation
(\ref{eqn:partfrac}).  Here we show how to reparametrize it so that
the resulting $2m$ parameters have a meaning of finite difference
coefficients.

A classical result of Stieltjes says that any rational function
of the form (\ref{eqn:partfrac}) with negative poles and positive
residues admits a representation as a Stieltjes continued fraction 
with positive coefficients
\begin{equation}
Y_m(s) = \frac{1}{\widehat{\kappa}_1 s + \dfrac{1}{\kappa_1 +
    \dfrac{1}{\widehat{\kappa}_1 s + \dfrac{1}{\ddots \;
        \dfrac{1}{\widehat{\kappa}_m s + \dfrac{1}{\kappa_m}}}}}}.
\label{eqn:cfrac}
\end{equation}
Obviously, this is true in our case, since $-\theta_j$ are the 
Ritz values of a negative definite operator $A$ and the residues
are given as squares in (\ref{eqn:thetac}).

Furthermore, it is known from \cite{druskin2000gaussian} that
(\ref{eqn:cfrac}) is a boundary response $w_1(s)$
(Neumann-to-Dirichlet map) of a second order finite difference scheme
with three point stencil
\begin{equation}
\frac{1}{\widehat{\kappa}_j} \left( \frac{w_{j+1}-w_j}{\kappa_j} -
\frac{w_j-w_{j-1}}{\kappa_{j-1}} \right) - s w_j = 0, \quad
j=2,\ldots,m, 
\label{eqn:fdw}
\end{equation}
and boundary conditions
\begin{equation}
\frac{1}{\widehat{\kappa}_1}\left( \frac{w_2-w_1}{\kappa_1} \right) -
s w_1 + \frac{1}{\widehat{\kappa}_1} = 0, \quad w_{m+1} = 0.
\label{eqn:bcw}
\end{equation}
These equations closely resemble the Laplace transform of 
(\ref{eqn:semidfwd}), except that the discretization is at $m$ 
nodes, which is much smaller than the dimension $N$ of the fine grid.

It is convenient to work henceforth with the logarithms $\log
\kappa_j$ and $\log \widehat{\kappa}_j$ of the finite difference
coefficients.  We now introduce the two mappings $\cQ$ and $\cR$ that
play the crucial role in our inversion method.  We refer to the first
mapping $\cQ : C(0, +\infty) \to \mathbb{R}^{2m}$ as \emph{data
  fitting}. It takes the time-dependent data $d(t)$ to the $2m$
logarithms of reduced order model parameters
\begin{equation}
\cQ(d(\,\cdot\,)) = \{ (\log \kappa_j, \log \widehat{\kappa}_j) \}_{j=1}^{m}
\end{equation}
via the following chain of mappings
\begin{equation}
\cQ \, : \, d(t) \stackrel{\mbox{\tiny (a)}}{\to} Y(s)
\stackrel{\mbox{\tiny (b)}}{\to} Y_m(s) \stackrel{\mbox{\tiny
    (c)}}{\to} \{ (c_j, \theta_j) \}_{j=1}^{m} \stackrel{\mbox{\tiny
    (d)}}{\to} \{ (\kappa_j, \widehat{\kappa}_j) \}_{j=1}^{m}
\stackrel{\mbox{\tiny (e)}}{\to} \{(\log \kappa_j, \log
\widehat{\kappa}_j) \}_{j=1}^{m}.
\label{eqn:datachain}
\end{equation}
Here step (a) is a Laplace transform of the measured data $d(t)$,
step (b) requires solving the rational interpolation problem
(\ref{eqn:matching}), which is converted to a partial fraction
form in step (c), which in turn is transformed to a Stieltjes
continued fraction form in step (d) with a variant of Lanczos
iteration, as explained in Appendix \ref{app:compjac}.

Note step (b) is the only ill-conditioned computation in the
chain. The ill-posedness is inherited from the instability of the
parabolic inverse problem and we mitigate it by limiting $m$. The
instability may be understood intuitively by noting that step (b) is
related to analytic continuation. We give more details in section
(\ref{subsec:ratinterp}).  In practice we choose $m$ so that the
resulting $\kappa_j$ and $\widehat{\kappa}_j$ are positive. This
determines the maximum number of degrees of freedom that we can
extract from the data at the present noise level.

The mapping $\cR : \mathbb{R}^N_+ \to \mathbb{R}^{2m}$  is 
the \emph{non-linear preconditioner}.
It takes the vector of
discrete resistivities $\br \in \mathbb{R}^N_+$ to the same output as 
$\cQ$. In simplest terms, it is a composition of $\cQ$ and the 
forward map (\ref{eqn:forwardmap}) $\cR = \cQ \circ \cF$ given by
\begin{equation}
\cR(\br) = \cQ(y(\,\cdot\,; \br)).
\end{equation}
Unlike the data fitting, the computation of $\cR$ can be done in a
stable manner using a chain of mappings
\begin{equation}
\cR: \; \br \stackrel{\mbox{\tiny (a)}}{\to} A(\br)
\stackrel{\mbox{\tiny (b)}}{\to} V \stackrel{\mbox{\tiny (c)}}{\to}
A_m \stackrel{\mbox{\tiny (d)}}{\to} \{ (c_j, \theta_j) \}_{j=1}^{m}
\stackrel{\mbox{\tiny (e)}}{\to} \{ (\kappa_j, 
\widehat{\kappa}_j) \}_{j=1}^{m}\stackrel{\mbox{\tiny (f)}}{\to} \{
(\log \kappa_j, \log \widehat{\kappa}_j) \}_{j=1}^{m}.
\label{eqn:kappachain}
\end{equation}
Here step (a) is just the definition of $A$, in (b) we compute the
orthonormal basis $V$ for the rational Krylov subspace
(\ref{eqn:ratkrylov}) on which $A$ is projected in step (c) to obtain
$A_m$ and $\bb_m$. Then, the poles and residues (\ref{eqn:thetac}) are
computed. The last two steps are the same as in the computation of
$\cQ$.

\subsection{Non-linearly preconditioned optimization}
\label{subsec:optim}
Now that we defined the data fitting $\cQ$ and the non-linear
preconditioner $\cR$ mappings, we can formulate our method for solving
the semi-discrete inverse problem as an optimization. Given the data
$d(t)$ for $t \in [0, +\infty)$ (recall (\ref{eqn:semiddata})),  we
  estimate the true resistivity $\br^{\mbox{\tiny true}}$ by
  $\br^\star$, the solution of the non-linear optimization problem
\begin{equation}
\br^\star = \mathop{\mbox{arg min}}\limits_{\br \in \mathbb{R}^N_+}
\frac{1}{2} \left\| \cQ(d(\,\cdot\,)) - \cR(\br))
\right\|_2^2.
\label{eqn:roptim}
\end{equation}

This is different than the typical optimization-based inversion, which
minimizes the $L_2$ norm of the (possibly weighted) misfit between the
measured data $d(\,\cdot\,)$ and the model $\cF(\br)$. Such an
approach is known to have many issues. In particular, the functional
is often non-convex with many local minima, which presents a challenge
for derivative-based methods (steepest descent, non-linear conjugate
gradients, Gauss-Newton). The convergence if often slow and some form
of regularization is required. In our approach we aim to convexify the
objective functional by constructing the non-linear preconditioner
$\cR$.

To explain the map $\cR$, let us give a physical interpretation of the
reduced model parameters $\{(\kappa_j,\widehat \kappa_j)\}_{j=1}^m$ by
introducing the change of coordinates $x \leadsto \xi$, so that
$r^{1/2} \partial_x =\partial_\xi$. The equation for $U(s,\xi)$, the
Laplace transform of $u(t,x(\xi))$, is
\begin{equation}
r^{-1/2} \partial_\xi \left( r^{1/2} \partial_\xi U\right) -s U = 0
\label{eqn:xi}
\end{equation}
and (\ref{eqn:fdw}) is its discretization on a staggered grid. The
coefficients $\widehat \kappa_j$ and $\kappa_j$ are the increments in the
discretization of the differential operators $r^{-1/2} \partial_\xi$
and $r^{1/2}\partial_\xi$, so they are proportional to the local
values of $r^{1/2}$ and $r^{-1/2}$, respectively. Since 
$\cR = \cQ \circ \cF$, an ideal choice of $\cQ$ would be an approximate 
inverse of $\cF$ for resistivity functions that can be approximated well 
on the discretization grid. It would interpolate in some manner the grid
values of $r(x)$, which are defined by $\{\widehat \kappa_j^2,
\kappa_j^{-2}\}_{j=1}^m$, up to some scaling factors that define the grid
spacing. 

Not all grids give good results, as explained in more detail in the
next section. We can factor out the unknown
grid spacings by working with the logarithm of the coefficients
$\{\widehat \kappa_j^2, \kappa_j^{-2}\}_{j=1}^m$ instead of the
resistivity $r$. Although it is possible to calculate a good grid, 
we do not require it in our inversion method. However, the grids
can be useful for judging the quality of different matching conditions
and for visualizing the behavior of $\cR$
as shown in section \ref{subsec:optgrid}.

The ill-posedness of the inverse problem is localized in the
computation of the data fitting term $\cQ(d(\,\cdot\,))$ that is
computed once.  The instability of this computation can be controlled
by reducing the size $m$ of the reduced order model projection
subspace.  A good strategy to follow in practice is to choose the
largest $m$ such that all the coefficients $\kappa_j$,
$\widehat{\kappa}_j$, $j=1,\ldots,m$ computed by $\cQ$ are positive.

Conventional optimization-based methods regularize the inversion by
adding a penalty term to the objective function. Our approach is
different. We acknowledge that there is a resolution vs. stability
trade-off by reducing the size of the reduced model, and view
regularization only as a means of adding prior information about the
unknown resistivity. If such information is available, we can
incorporate it at each iteration of the Gauss-Newton method via a
correction in the null space of the Jacobian $\cD \cR$. The
regularization scheme is discussed in detail in section
\ref{subsec:regular}.

Inversion via model order reduction is a recent idea that was
introduced in \cite{dsz2012,druskin2007combining,beattietal}.  In
particular, the approach in \cite{dsz2012} uses maps
$\cQ_{_{c\theta}}$ and $\cR_{_{c\theta}}$ to the spectral parameters
of the reduced order model $\theta_j$ and $c_j$, $j=1,\ldots,m$.
Unlike the continued fraction coefficients $\kappa_j$ and
$\widehat{\kappa}_j$, the poles $-\theta_j$ and residues $c_j$ do not
have a physical meaning of resistivity and the mapping
$\cR_{_{c\theta}}$ does not behave like an approximate identity, as
$\cR$ does. Our definition of the mapping $\cR$ is based on the ideas
from \cite{borcea2005continuum}. It allows
us to improve the results of \cite{dsz2012}. The improvement
becomes especially pronounced in cases of high resistivity contrast,
as shown in the numerical comparison in section \ref{sec:numres}.


The idea of convexification of the non-linear inverse problem has been
pursued before for hyperbolic equations in 
\cite{beretta2014inverse, klibanov2007numerical} and global
convergence results were obtained in \cite{beilina2008globally}. 
However, it is not clear if these approaches apply to parabolic equations. 
Although both hyperbolic and parabolic equations can be transformed to 
the frequency domain via Fourier and Laplace transforms, their data is 
mapped to different parts of the complex plane where the spectral 
parameter lies. It is known that the transformations from the hyperbolic 
to the parabolic data are unstable \cite{kabanikhin2011inverse}, so 
one cannot directly apply the methods from 
\cite{beilina2008globally, beretta2014inverse, klibanov2007numerical}
to the parabolic problem. The model reduction approach in this paper 
gives a specially designed discrete problem of small size which can be 
inverted stably.

\subsection{Connection to optimal grids}
\label{subsec:optgrid}

We explain here that the map $\cR$ relates to an interpolant of the
resistivity $r$ on a special, so-called optimal grid.  Although our
inversion method does not involve directly such grids, it is
beneficial to study them because they provide insight into the
behavior of the non-linear preconditioner $\cR$.  We give a brief
description of the grids, and show that they depend strongly on the
interpolation nodes $s_j$ in the matching conditions
(\ref{eqn:matching}). We use this dependence to conclude that not all
rational approximants of the transfer function $Y(s;\br)$ are useful
in inversion, as shown in section \ref{subsec:matchcond}.  Then, we
use in section \ref{subsec:precondaction} the optimal grids and the
continued fraction coefficients $\kappa_j$ and $\widehat{\kappa}_j$ to
visualize the action of $\cR$ on $\br$. This allows us to display how the
non-linear preconditioner approximates the identity.

\subsubsection{Optimal grids}

The optimal grids have been introduced in \cite{druskin2000gaussian, 
druskin2000gsr, IngDruKni} to obtain
exponential convergence of approximations of the Dirichlet-to-Neumann
map. They were used and analyzed in the context of inverse spectral
problems in \cite{borcea2005continuum} and in the inverse problem of
electrical impedance tomography in \cite{BorDruGue, BDMG-10}. 
The grids are defined by the coefficients of the
continued fraction representation of the rational response function
corresponding to a reduced model of a medium with constant resistivity
$r^{(0)}$.  In principle we can define grids for other reference
resistivities, not necessarily constant, but the results in 
\cite{BorDruGue, borcea2005continuum, BDMG-10}
show that the grids change very little with respect to $r(x)$.  This
is why they are very useful in inversion.

Let then $r^{ (0)}(x) \equiv 1$ so that the change of coordinates
in (\ref{eqn:xi}) is trivial ($\xi = x$) and obtain from
(\ref{eqn:fdw}-\ref{eqn:bcw}) that $ \{(\kappa_j^{ (0)},
\widehat{\kappa}_j^{ (0)})\}_{j=1}^{m}$ are the optimal grid
steps in a finite difference scheme for the equation
\[
\frac{\partial^2 w}{\partial x^2} - s w = 0.
\]  
The primary and dual grid points are
\begin{equation}
x^{ (0)}_j = \sum_{k=1}^{j} \kappa_k^{ (0)}, \quad
\widehat{x}^{ (0)}_j = \sum_{k=1}^{j} \widehat{\kappa}_k^{
  (0)}, \quad j=1,\ldots,m,
\label{eqn:xopt}
\end{equation}
with boundary nodes $x^{ (0)}_0 = \widehat{x}^{ (0)}_0 = 0$.
In the numerical experiments we observe that the grid is staggered,
i.e. the \emph{primary} nodes $x^{ (0)}_j$ and the \emph{dual}
nodes $\widehat{x}^{ (0)}_j$ obey the interlacing conditions
\begin{equation}
0 = \widehat{x}^{ (0)}_0 = x^{ (0)}_0 < 
\widehat{x}^{ (0)}_1 < x^{ (0)}_1 <
\widehat{x}^{ (0)}_2 < x^{ (0)}_2 < \ldots <
 x^{ (0)}_{m-1}<
\widehat{x}^{ (0)}_m < x^{ (0)}_m \leq 1.
\label{eqn:xinterlace}
\end{equation}
We do not prove (\ref{eqn:xinterlace}) here, although it is possible
to do so at least in some settings. 

The optimal grids are closely connected with the \emph{sensitivity
 functions}, which for the semi-discrete problem are given by the rows
of the Jacobian $\cD \cR \in \mathbb{R}^{2m \times N}$ defined as
\begin{equation}
\left( \cD \cR \right)_{j,k} = \left\{ \begin{tabular}{ll}
  $\dfrac{\partial \log \kappa_j}{\partial r_k}$, & if $1 \leq j \leq m$
  \\ \\ $\dfrac{\partial \log \widehat{\kappa}_j}{\partial r_k}$, & if $m+1
  \leq j \leq 2m$
\end{tabular}
\right., \quad k=1,\ldots,N.
\nonumber
\end{equation}
The studies in \cite{BorDruGue, BDMG-10} show that the sensitivity function corresponding to 
$\kappa_j$ is localized around the corresponding grid cell 
$(\widehat{x}^{(0)}_j, \widehat{x}^{(0)}_{j+1})$, and its maximum 
is near $x^{(0)}_j$.  The same holds for $\widehat{\kappa}_j$, after
interchanging the primary and dual grid nodes. Moreover, the columns
of the pseudoinverse $(\cD \cR)^\dagger$ have similar localization
behavior.  The Gauss-Newton update is the linear combination of the
columns of $(\cD \cR)^\dagger$, and therefore optimal grids are useful
for inversion. They localize well features of the resistivity that are
recoverable from the measurements.

A good grid should have two properties. First, it should be refined
near the point of measurement to capture correctly the loss of
resolution away from $\mathcal{B}_A$. Second, the nodes should not 
all be clustered near $\mathcal{B}_A$ because when the nodes get 
too close, the corresponding rows of $\cD \cR$ become almost linearly
dependent, and the Jacobian is ill-conditioned.  

\subsubsection{Matching conditions}
\label{subsec:matchcond}

We study here the grids for three choices of matching conditions
(\ref{eqn:matching}). The first corresponds to $l=1$, $s_1=0$ and
$M_1=m$ (simple Pad\'{e}) and yields the rational Krylov subspace
\[
\mathcal K_m(0) = \mbox{span} \left\{ A^{-1} \bb, A^{-2} \bb, \ldots,
A^{-m} \bb \right\}.
\]
This approximant has the best accuracy near the interpolation point
($s=0$), and is obviously inferior for global approximation in $s$
when compared to multipoint Pad\'{e} aproximants.  The other two
choices match $Y(s;\br)$ and its first derivative at nodes
$\widetilde{\bf s} = (\widetilde{s}_j)_{j=1}^m$ that are distributed
geometrically
\begin{equation}
\widetilde{s}_j = \widetilde{s}_1 
\left( \frac{\widetilde{s}_2}{\widetilde{s}_1} \right)^{j-1}, 
\quad j = 1,\ldots,m.
\label{eqn:geometric}
\end{equation}
We use henceforth the tilde to distinguish these nodes from those
obtained with the change of variables
\begin{equation}s_j = \frac{\widetilde{s}_j }{
\widetilde{s}_m} \in (0,1],
\label{eqn:changesvar}
\end{equation}
intended to improve the conditioning of the interpolation.  The
mathching conditions at $\widetilde{\bf s}$ yield the rational Krylov
projection subspace
\[
\mathcal K_m(\widetilde{\bf s}) = 
\mbox{span} \left\{  
(\widetilde{s}_1 I - A)^{-1} \bb, \ldots,
(\widetilde{s}_m I - A)^{-1} \bb \right\},
\]
and the two choices of interpolants differ in the rate of growth of
$\widetilde{\bf s}$.  The second interpolant uses the rate of growth
$\widetilde s_2/\widetilde s_1 = 1+12/m$ in (\ref{eqn:geometric}) and
$\widetilde{s}_1 = 2$, so that
\begin{equation}
\widetilde{s}_j = 2 \left( 1 + \frac{12}{m}
\right)^{j-1}, \quad j = 1, \ldots, m.
\label{eqn:sigmatilde}
\end{equation}
This choice approximates the Zolotarev nodes \cite{IngDruKni} which
arise in the optimal rational approximation of the transfer function
$Y(s;\br)$ over a real positive and bounded interval of $s$.  The
third interpolant uses a faster rate of growth of $\widetilde {\bf s}$
and gives worse results, as illustrated in the numerical experiment
below.

\begin{figure}[t!]
\begin{center}
\includegraphics[width=0.9\textwidth]{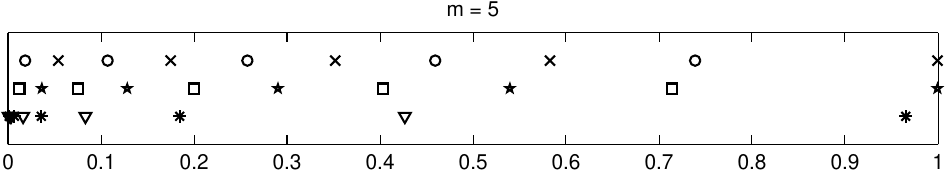} \vskip10pt
\includegraphics[width=0.9\textwidth]{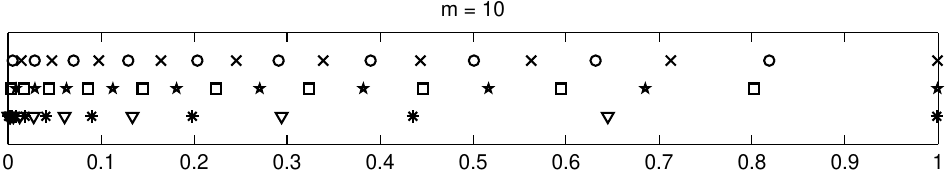}
\caption{Primary $x^{ (0)}_j$ and dual $\widehat{x}^{
    (0)}_j$ optimal grid nodes for $j=1,\ldots,m$ ($m=5,10$) and
  different choices of matching conditions.  Moment matching at $s=0$:
  primary $\times$, dual $\circ$.  Interpolation at geometrically
  distributed interpolation nodes: 
  primary $\star$ and dual $\square$ for slowly growing $\widetilde{\bf s}$;
  primary $*$ and dual $\bigtriangledown$ for fast growing $\widetilde{\bf s}$.  
  The number of fine grid steps in the semi-discretized model is $N=1999$.}
\label{fig:grids}
\end{center}
\end{figure}

\begin{figure}[t!]
\begin{center}
\includegraphics[width=0.55\textwidth]{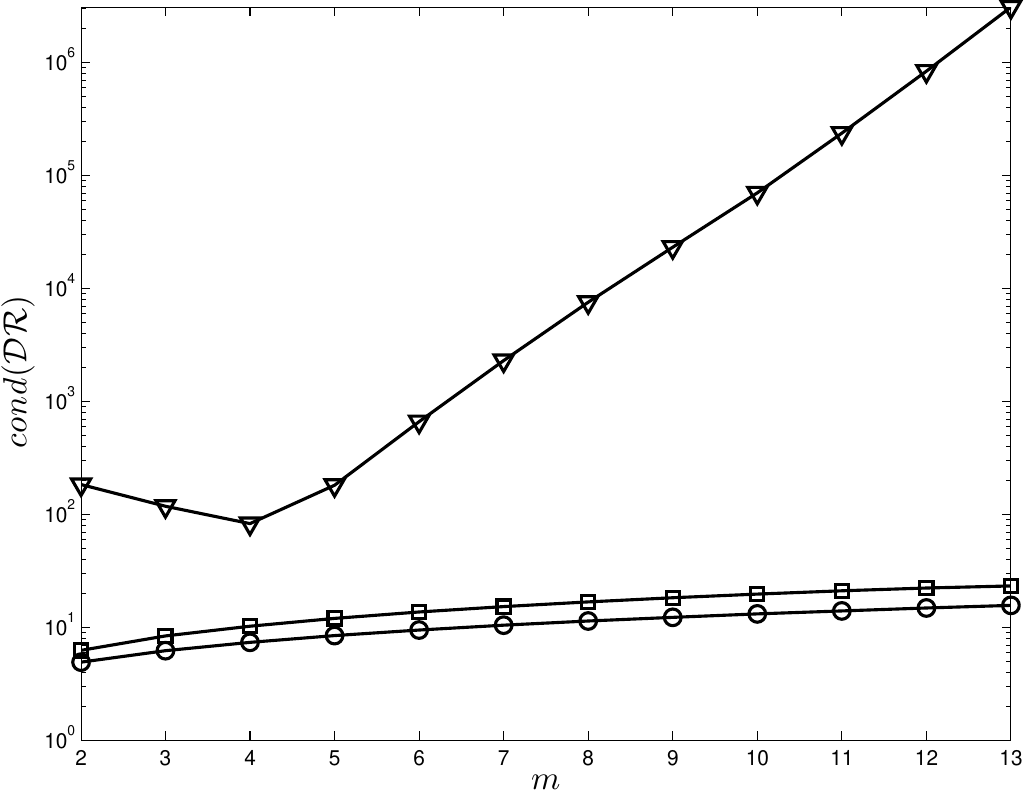}
\caption{Dependence of the
  condition number of $\cD \cR$ on the size $m$ of the reduced model
  for different matching conditions. Moment matching at $s=0$ is
  $\circ$, interpolation at geometrically distributed nodes: slowly
  growing $\widetilde{\bf s}$ is $\square$, fast growing
  $\widetilde{\bf s}$ is $\bigtriangledown$.  The number of fine grid
  steps in the semi-discretized model is $N=1999$. }
\label{fig:condJ}
\end{center}
\end{figure}

We show the optimal grids for all three choices of matching conditions
in Figure \ref{fig:grids} for reduced models of sizes $m=5,10$. We
observe that the nodes of the grids corresponding to fast growing
$\widetilde{\bf s}$ are clustered too close to the measurement point
$x=0$. Thus, inversion results are expected to have poor resolution
throughout the rest of the domain away from the origin. In addition,
the clustering of the grid nodes leads to poor conditioning of the
Jacobian $ \cD \cR$.  This is illustrated in Figure \ref{fig:condJ},
where the condition numbers are plotted against the size $m$ of the
reduced model. We observe that the condition number of the Jacobian
for the reduced model with fast growing $\widetilde{\bf s}$ increases
exponentially.  The condition numbers of the Jacobians for the other
two choices of matching conditions grow very slowly.

We can explain why the case of fast growing $\widetilde{\bf s}$ is
undesirable in inversion by looking at the limiting case of
approximation at infinity. The simple Pad\'{e} approximant at infinity
corresponds to the Krylov subspace
\[
\cK_m(+\infty) = \mbox{span} 
\left\{ \bb, A \bb, \ldots, A^{m-1} \bb \right\}
\]
which is unsuitable for inversion in our setting. To see this, recall from 
(\ref{eqn:Ar}) that $A(\br)$ is tridiagonal and $\bb$ is a scalar multiple 
${\bf e}_1$. Thus, for any $j \in \mathbb{Z}_+$ only the first $j+1$ 
components of the vector $A^j \bb$ are non-zero. When $A(\br)$ is 
projected on $K_m(+\infty)$ in (\ref{eqn:rom-project}), the reduced 
model matrix $A_m$ is aware only of the upper left $(m+1)\times(m+1)$ block 
of $A(\br)$, which depends only on the first $m+1$ entries of $\br$. 
This is unsuitable for inversion where we want $A_m$ to capture the 
behavior of the resistivity in the whole interval, even for small $m$.  
The corresponding optimal grid steps will simply coincide with the first
$m$ grid steps $h$ of the fine grid discretization in (\ref{eqn:Ar}).

When the interpolation nodes grow too quickly, we are near the
limiting case of simple Pad\'{e} approximant at infinity, and 
the first optimal grid steps are 
$
\widehat{x}^{ (0)}_1 \approx x^{ (0)}_1 \approx h.
$
Consequently, the rows 
$\left( \dfrac{\partial \kappa_1}{\partial r_k} \right)_{k=1}^{N}$ and 
$\left( \dfrac{\partial \widehat{\kappa}_1}{\partial r_k} \right)_{k=1}^{N}$
of $\cD \cR$ are almost collinear and the Jacobian is poorly conditioned, as
shown in Figure \ref{fig:condJ}. 

\subsubsection{Action of the non-linear preconditioner on the resistivity}
\label{subsec:precondaction}

The optimal grids can be used to obtain resistivity reconstructions
directly, without optimization, as was done in
\cite{borcea2005continuum} for the inverse spectral problem and in
\cite{BorDruGue, BDMG-10} for electrical
impedance tomography.  We do not use this approach, but we show here
such reconstructions to display the behavior of the non-linear
preconditioner $\cR$ when acting on $\br$.

\begin{figure}[t!]
\begin{center}
\begin{tabular}{ccc}
$r_Q$ & $r_L$ & $r_J$ \\
\hskip-0.01\textwidth
\includegraphics[width=0.32\textwidth]{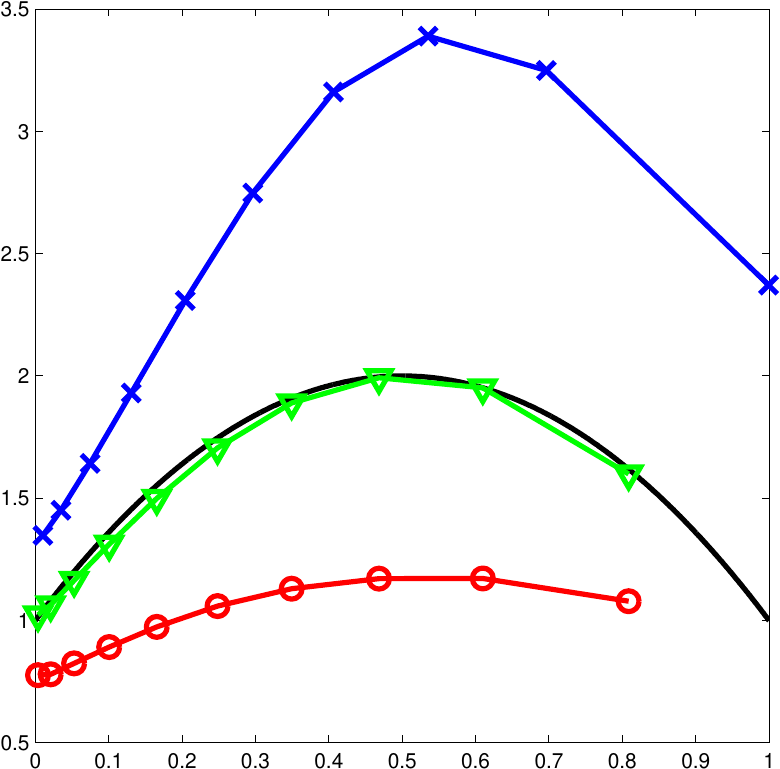}
&
\hskip-0.01\textwidth
\includegraphics[width=0.32\textwidth]{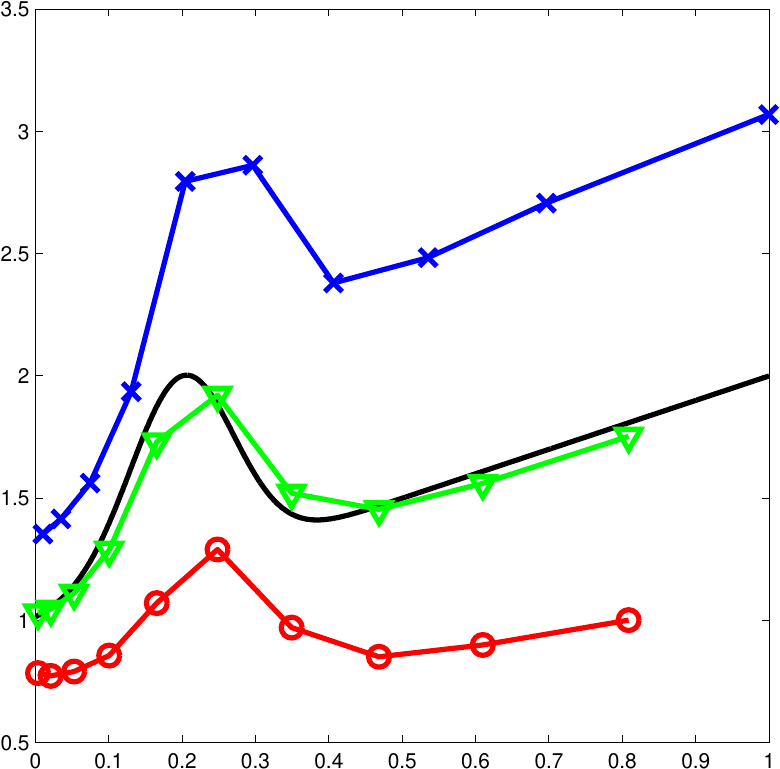}
&
\hskip-0.01\textwidth
\includegraphics[width=0.32\textwidth]{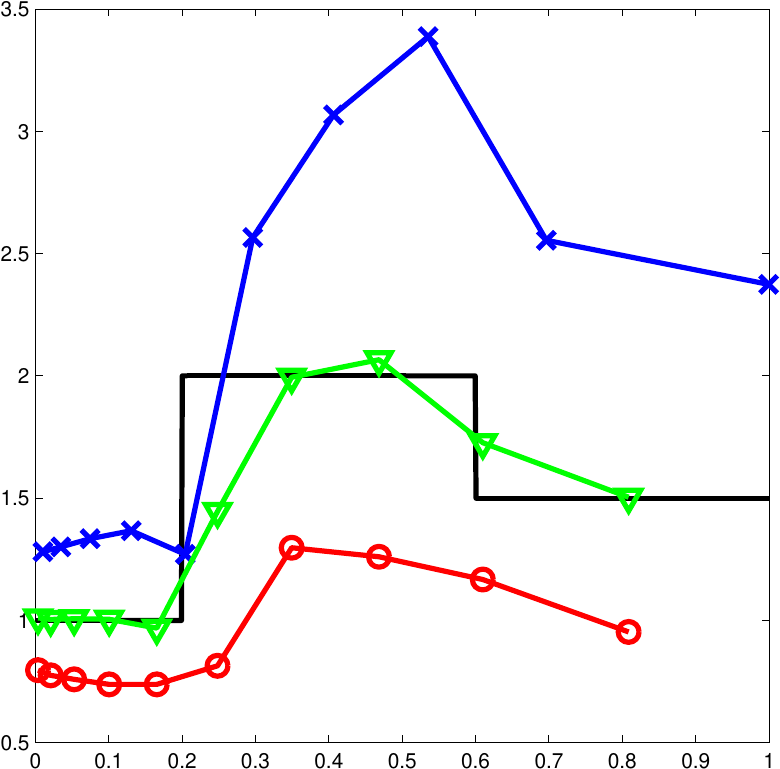}
\end{tabular}
\caption{Action of the non-linear preconditioner $\cR$ on the
  resistivities $r_Q$, $r_L$ and $r_J$ (solid black line) defined in
  (\ref{eqn:rQL}) and (\ref{eqn:rJ}). The ``primary'' ratios
  $(x^{(0)}_j, \zeta_j)$ are blue $\times$, the ``dual'' ratios
  $(\widehat{x}^{(0)}_j, \widehat{\zeta}_j)$ are red $\circ$ for
  $j=1,\ldots,m$ with $m=10$. The geometric averages
  $(\widehat{x}^{(0)}_j, \widetilde{\zeta}_j)$ are green
  $\bigtriangledown$.  }
\label{fig:precond}
\end{center}
\end{figure}

\begin{figure}[ht!]
\begin{center}
\includegraphics[width=0.5\textwidth]{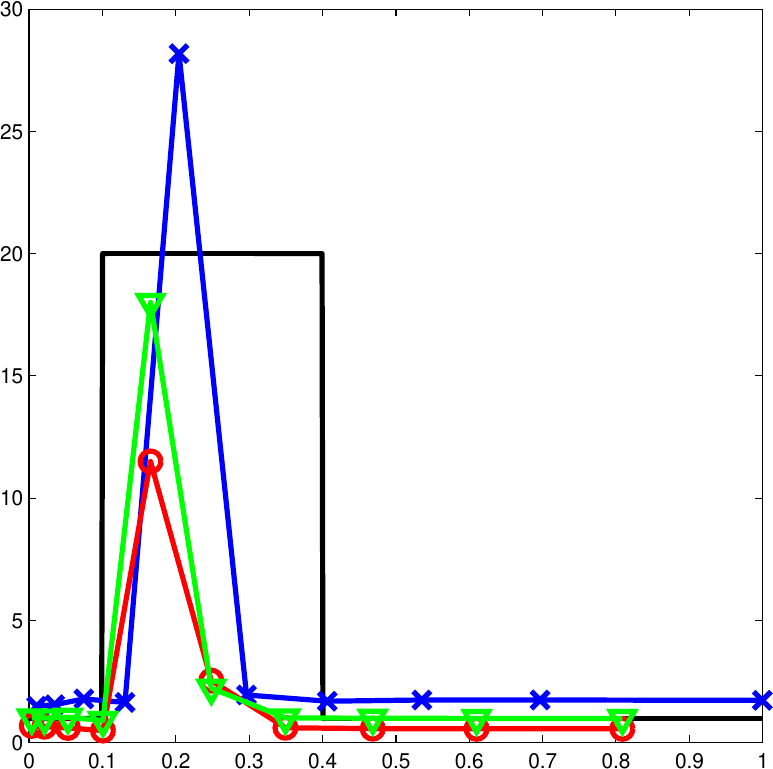}
\caption{Action of the non-linear preconditioner $\cR$ on a piecewise
  constant resistivity of contrast $20$. Same setting as in figure
  \ref{fig:precond}.  }
\label{fig:prehigh}
\end{center}
\end{figure}

Recall from  equation (\ref{eqn:xi}) and the explanation in 
section \ref{subsec:optim} that $\widehat{\kappa}_j$ and $\kappa_j$ 
are proportional to the values of $r^{1/2}(x)$ and $r^{-1/2}(x)$ 
around the corresponding optimal grid points. If we take as the
proportionality coefficients the values $\widehat{\kappa}_j^{(0)}$
and $\kappa_j^{(0)}$ then we expect the ratios
\begin{equation}
\zeta_j = \left( \frac{\kappa_j^{(0)}}{\kappa_j} \right)^2, \quad
\widehat{\zeta}_j = \left(
\frac{\widehat{\kappa}_j}{\widehat{\kappa}_j^{(0)}} \right)^2, \quad
j=1,\ldots,m,
\label{eqn:zeta}
\end{equation}
to behave roughly as $r(x_j^{(0)})$ and $r(\widehat{x}_j^{(0)})$.  In
practice, a more accurate estimate of the resistivity can be obtained
by taking the geometric average of $\zeta_j$ and $\widehat{\zeta}_j$
\begin{equation}
\widetilde{\zeta}_j = \sqrt{\zeta_j \widehat{\zeta}_j} = \frac{
  \kappa_j^{(0)} \widehat{\kappa}_j }{ \kappa_j
  \widehat{\kappa}_j^{(0)} }, \quad j=1,\ldots,m.
\label{eqn:zetatilde}
\end{equation}
Since building a direct inversion algorithm is not our focus, we only
show (\ref{eqn:zetatilde}) for comparison purposes.

In figure \ref{fig:precond} we display the ratios (\ref{eqn:zeta})
plotted at the nodes of the optimal grid $(x^{(0)}_j, \zeta_j)$,
$(\widehat{x}^{(0)}_j, \widehat{\zeta}_j)$, $j=1,\ldots,m$. We take
the same resistivities $r_Q$, $r_L$ and $r_J$ that are used in the
numerical experiments in section \ref{sec:numres}. They are defined in
(\ref{eqn:rQL}) and (\ref{eqn:rJ}).  We observe that the curve defined
by the linear interpolation of the ``primary'' points $(x^{(0)}_j,
\zeta_j)$ overestimates the true resistivity, while the ``dual'' curve
passing through $(\widehat{x}^{(0)}_j, \widehat{\zeta}_j)$
underestimates it.  Both curves capture the shape of the resistivity
quite well, so when taking the geometric average (\ref{eqn:zetatilde})
the reconstruction falls right on top of the true resistivity. This
confirms that $\cR$ resolves most of the non-linearity of the problem
and thus acts on the resistivity as an approximate identity.

We can also illustrate how well $\cR$ resolves the non-linearity of
the problem by considering an example of high contrast resistivity. In
figure \ref{fig:prehigh} we plot the same quantities as in figure
\ref{fig:precond} in the case of piecewise constant resistivity of
contrast $20$. The contrast is captured quite well by
$\widetilde{\zeta}_j$, while the shape of the inclusion is shrunk.
This is one of the reasons why we use $\cR$ as a preconditioner for
optimization and not a reconstruction mapping. Optimization allows us
to recover resistivity features on scales that are smaller than those
captured in $\zeta_j$ and $\widehat{\zeta}_j$. Moreover, since $\cR$
resolves most of the non-linearity of the inverse problem, the
optimization avoids the pitfalls of traditional data
misfit minimization approaches, such as sensitivity to the initial
guess, numerous local minima and slow convergence.

\subsection{Data fitting via the rational interpolation}
\label{subsec:ratinterp}

Unlike $\cR(\br) = \cQ(y(\,\cdot\,;\br))$ computed using the chain of
mappings (\ref{eqn:kappachain}) with all stable steps, the computation
of the data fitting term $\cQ(d(\,\cdot\,))$ requires solving an
osculatory rational interpolation problem (\ref{eqn:matching}) in step
(b) of (\ref{eqn:datachain}) to obtain the rational interpolant
$Y_m(s)$ of the transfer function $Y(s)$.  This involves the solution
of a linear system of equations with an ill-conditioned matrix, or
computing the singular value decomposition of such matrix. We use the
condition number of the matrix to assess the instability of the
problem. The condition number grows exponentially with $m$, but the
rate of growth depends on the matching conditions used in the rational
interpolation.  We show this for the two choices of matching
conditions: interpolation of $Y(s)$ and its first derivatives at
distinct nodes $\widetilde{\bf s}$ distributed as in
(\ref{eqn:sigmatilde}) (multipoint Pad\'{e}), and matching of moments
of $Y(s)$ at $s=0$ (simple Pad\'{e}).  We describe first both Pad\'{e}
interpolation schemes and then compare their stability numerically.

\vspace{0.1in} \textbf{Multipoint Pad\'{e}:} Let us rewrite the
reduced order transfer function (\ref{eqn:Gmtransfer}) in the form
\begin{equation}
Y_m(s) = \frac{f(s)}{g(s)} = \frac{f_0 + f_1 s + \ldots + f_{m-1}
  s^{m-1}}{g_0 + g_1 s + \ldots + g_{m} s^{m}},
\label{eqn:ratio}
\end{equation}
where $f(s)$ and $g(s)$ are polynomials defined up to a common
factor. We use this redundancy later to choose a unique solution of an
underdetermined problem.  The matching conditions (\ref{eqn:matching})
of $Y(s)$ and $Y'(s)$ at the distinct interpolation nodes $0 <
\widetilde{s}_1 < \widetilde{s}_2 < \ldots < \widetilde{s}_m$ are
\begin{equation}
\left\{ \begin{tabular}{l} $f(\widetilde{s}_j) -
  Y_m(\widetilde{s}_j) g(\widetilde{s}_j) = 0$
  \\ $f'(\widetilde{s}_j) - Y_m'(\widetilde{s}_j)
  g(\widetilde{s}_j) - Y_m(\widetilde{s}_j)
  g'(\widetilde{s}_j) = 0$
\end{tabular} \right., \quad j=1,\ldots,m.
\label{eqn:osculatory}
\end{equation}
Next, we recall the change of variables (\ref{eqn:changesvar}) and
define the Vandermonde-like $m \times (m+1)$ matrices
\begin{equation}
\mathcal{S} = \begin{bmatrix} 1 & s_1 & s_1^2 & \ldots &
  s_1^m \\ 1 & s_2 & s_2^2 & \ldots & s_2^m
  \\ \vdots & \vdots & \vdots & \vdots & \vdots \\ 1 & s_m &
  s_m^2 & \ldots & s_m^m
\end{bmatrix}, \quad
\mathcal{S}' = \frac{1}{\widetilde{s}_m } \begin{bmatrix} 
0 & 1 & 2s_1 & \ldots & ms_1^{m-1} \\
0 & 1 & 2s_2 & \ldots & ms_2^{m-1} \\
\vdots & \vdots & \vdots & \vdots & \vdots \\
0 & 1 & 2s_m &  \ldots & ms_m^{m-1}
\end{bmatrix},
\end{equation}
and the diagonal matrices
\begin{equation}
\mathcal{Y} = \mbox{diag} \left( Y(\widetilde{s}_1), \ldots,
Y(\widetilde{s}_m) \right), \quad \mathcal{Y}' = \mbox{diag} \left(
Y'(\widetilde{s}_1), \ldots, Y'(\widetilde{s}_m) \right).
\label{eqn:FFprime}
\end{equation}
This allows us to write equations (\ref{eqn:osculatory}) in
matrix-vector form as an underdetermined problem
\begin{equation}
\mathcal{P} {\bf u} = 0, \quad {\bf u} \in \mathbb{R}^{2m+1}, 
\label{eqn:rv}
\end{equation}
with 
\begin{equation}
\label{eqn:Rm}
\mathcal{P} = \begin{bmatrix} \mathcal{S}_{1:m,\,1:m} & - \mathcal{Y}
  \mathcal{S} \\ \mathcal{S}_{1:m,\,1:m}' & - \mathcal{Y}' \mathcal{S} -
  \mathcal{Y} \mathcal{S}'
\end{bmatrix}
\in \mathbb{R}^{2m \times (2m+1)},
\end{equation}
and  
\begin{align}
f_j & =  \widetilde{s}_m^{-j} u_{j+1}, \quad j = 0,\ldots,m-1, \nonumber 
\\ g_j & =  \widetilde{s}_m^{-j} u_{j+m+1}, \quad j =
0,\ldots,m. \label{eqn:Rm1} 
\end{align}
The problem is underdetermined
because of the redundancy in (\ref{eqn:ratio}). We eliminate it by the
additional condition $\| {\bf u} \|_2 = 1$, which makes it possible to
solve (\ref{eqn:rv}) via the singular value decomposition. If we let 
$U$ be the matrix of right singular vectors of $\mathcal{P}$, then 
\begin{equation}
{\bf u} = U_{1:(2m+1),\,2m+1}.
\end{equation}

Once the polynomials $f$ and $g$ are determined from (\ref{eqn:Rm1}),
we can compute the partial fraction expansion
(\ref{eqn:partfrac}). The poles $-\theta_j$ are the roots of $g(s)$,
and the residues $c_j$ are given by
\begin{equation}
c_j = \frac{f(-\theta_j)}{g_m
  \mathop{\prod\limits_{k=1}^{m}}\limits_{k \neq j} (\theta_k -
  \theta_j)}, \quad j = 1, \ldots, m, \nonumber
\end{equation}
assuming that $\theta_j$ are distinct. Finally, $\kappa_j$ and
$\widehat{\kappa}_j$ are obtained from $\theta_j$ and $c_j$ via
a Lanczos iteration, as explained in Appendix \ref{app:compjac}.

\vspace{0.1in}
\textbf{Simple Pad\'{e}:} When $l=1$, $s_1 = 0$ and $M_1 = m$, we have
a simple Pad\'{e} approximant which matches the first $2m$ moments of
$y(t;\br)$ in the time domain, because
\begin{equation}
\left. \frac{\partial^j Y_m}{\partial s^j} \right|_{s=0} = 
\left. \frac{\partial^j Y}{\partial s^j} \right|_{s=0} = 
(-1)^j \int_{0}^{+\infty} y(t;\br) t^j dt, \quad j=0,1,\ldots,2m-1.
\nonumber
\end{equation}
A robust algorithm for simple Pad\'{e} approximation is proposed in
\cite{gonnet2011robust}. It is also based on the singular value
decomposition.  If $Y(s)$ has the Taylor expansion at $s = 0$
\begin{equation}
Y(s) = \tau_0 + \tau_1 s + \tau_2 s^2 + \ldots + \tau_{2m-1} s^{2m-1}
+ \ldots,
\label{eqn:Gtaylor}
\end{equation}
then the algorithm in \cite{gonnet2011robust} performs a singular value
decomposition of the  Toeplitz matrix
\begin{equation}
\mathcal{T} = \begin{bmatrix}
\tau_m      & \tau_{m-1}  & \cdots & \tau_1 & \tau_0 \\
\tau_{m+1}  & \tau_{m}    & \cdots & \tau_2 & \tau_1 \\
\vdots      & \vdots      & \ddots & \vdots & \vdots \\
\tau_{2m-1} & \tau_{2m-2} & \cdots & \tau_m & \tau_{m-1}
\end{bmatrix} \in \mathbb{R}^{m \times (m+1)}.
\end{equation}
If $U\in \mathbb{R}^{(m+1) \times (m+1)}$ is the matrix of right
singular vectors of $\mathcal{T},$ then the coefficients in
(\ref{eqn:ratio}) satisfy
\begin{equation}
\begin{pmatrix}
g_0 \\ g_1 \\ \vdots \\ g_m
\end{pmatrix} = U_{1:(m+1),\,m+1},
\end{equation}
and 
\begin{equation}
\begin{pmatrix}
f_0 \\ f_1 \\ \vdots \\ f_{m-1}
\end{pmatrix} = 
\begin{bmatrix}
\tau_0     & 0          & \cdots & 0      & 0      \\
\tau_1     & \tau_0     & \cdots & 0      & 0      \\
\vdots     & \vdots     & \ddots & \vdots & \vdots \\
\tau_{m-1} & \tau_{m-2} & \cdots & \tau_0 & 0      \\
\end{bmatrix}
\begin{pmatrix}
g_0 \\ g_1 \\ \vdots \\ g_m
\end{pmatrix}.
\label{eqn:pqpade}
\end{equation}
We refer the reader to \cite{gonnet2011robust} for detailed explanations.

\vspace{0.1in} \textbf{Comparison:} To compare the performance of the
two interpolation procedures, we give in Table \ref{tab:ratcond} the
condition numbers of matrices $\mathcal{P}$ and $\mathcal{T}$. They
are computed for the reference resistivity $\br^{(0)}$, for
several values of the size $m$ of the reduced model.  We observe that
while both condition numbers grow exponentially, the growth rate is
slower for the multipoint Pad\'{e} approximant. Thus, we conclude that
it is the best of the choices of matching conditions considered in
this section. It allows a more stable computation of $\cQ(d(\,\cdot\,))$, it
gives a good distribution of the optimal grid points and a
well-conditioned Jacobian $\cD \cR$.
\begin{table}[ht!]
\begin{center}
\caption{Condition numbers of $\mathcal{P}$ (multipoint Pad\'{e} approximant) and
  $\mathcal{T}$ (simple Pad\'{e} approximant). The number of fine grid steps
  in the semi-discretized model is $N=299$. }
\begin{tabular}{c|c|c|c|c|c}
$m$ & 2 & 3 & 4 & 5 & 6 \\ \hline cond($\mathcal{P}$) & $ 4.43 \cdot
10^{2}$ & $ 6.73 \cdot 10^{4}$ & $ 1.85 \cdot 10^{7}$ & $ 6.95 \cdot
10^{9}$ & $ 3.83 \cdot 10^{12}$ \\ \hline cond($\mathcal{T}$) & $ 5.28 \cdot
10^{1}$ & $ 1.26 \cdot 10^{5}$ & $ 1.84 \cdot 10^{9}$ & $ 9.14 \cdot
10^{13}$ & $ 2.86 \cdot 10^{16}$
\end{tabular}
\label{tab:ratcond}
\end{center}
\end{table}

\subsection{Regularization of Gauss-Newton iteration}
\label{subsec:regular}

Our method of solving the optimization problem (\ref{eqn:roptim}) uses
a Gauss-Newton iteration with regularization similar to that in
\cite{BorDruGue}. We outline it below, and refer to the next section
for the precise formulation of the inversion algorithm.


Recall that $\cR$ maps the vectors $\br \in \mathbb{R}^N_+$ of 
resistivity values on the fine grid to $2m$ reduced order model 
parameters $\{( \log \kappa_j, \log \widehat \kappa_j)\}_{j=1}^m$. 
Thus, the Jacobian $\cD \cR$ has dimensions $\mathbb{R}^{2m \times N}$. 
Since the reduced order model is much coarser than the fine grid 
discretization $2m \ll N$, the Jacobian $\cD \cR(\br)$ has a large 
null space. At each iteration the Gauss-Newton update to
$\br$ is in the $2m$ dimensional range of the pseudoinverse 
$\cD \cR^\dagger(\br)$ of the Jacobian. This leads to low
resolution in practice, because $m$ is kept small to mitigate the
sensitivity of the inverse problem to noise. If we have prior
information about the unknown $\br^{\mbox{\tiny true}}$, we can use it
to improve its estimate $\br^\star$. 

We incorporate the prior information about the true resistivity
into a penalty functional $\mathcal{L}(\br)$. For example, 
$\mathcal{L}(\br)$ may be the total variation norm of $\br$ if 
$\br^{\mbox{\tiny true}}$ is known to be piecewise constant, 
or the square of the $\ell_2$ norm of $\br$ or of its 
derivative if $\br^{\mbox{\tiny true}}$ is expected to be smooth.

In our inversion method we separate the minimization of the norm of
the residual $\cQ(d(\,\cdot\,))-\cR(\br)$ and the penalty functional
$\mathcal{L}(\br)$. At each iteration we compute the standard 
Gauss-Newton solution $\br$ and then we add a correction to obtain
a regularized iterate $\brho$. The correction $\brho-\br$ is in the 
null space of $\cD \cR$, so that the residual remains unchanged.  
We define it as the minimizer of the constrained optimization problem
\begin{equation}
\mathop{\mbox{minimize}}\limits_{\mbox{s.t. } [\cD \cR]
    (\br - \brho) = 0} \mathcal{L}(\brho),
\label{eqn:minreg}
\end{equation}
which we can compute explicitly in the case of a weighted discrete $H^1$
seminorm regularization, assumed henceforth, 
\begin{equation}
\mathcal{L}(\br) = \frac{1}{2} \| {W}^{1/2} \widetilde{D} \br \|_2^2.
\label{eqn:h1w}
\end{equation}
Here the matrix $\widetilde{D}$ is a truncation of $D$ defined as
\[
\widetilde{D} = D_{1:(N-1),\,1:N} \in \mathbb{R}^{(N-1) \times N},
\] and ${W} \in \mathbb{R}^{(N-1) \times (N-1)}$ is a diagonal matrix 
of weights.  We specify it below depending on the prior information on 
the true resistivity. 

With the choice of penalty in the form (\ref{eqn:h1w}) the optimization 
problem (\ref{eqn:minreg}) is quadratic with linear constraints, and 
thus $\brho$ can be calculated from the first order optimality 
conditions given by the linear system
\begin{align}
\widetilde{D}^T W \widetilde{D} \,  \brho + [\cD \cR]^T
\bflambda & =  0, \\ \relax[\cD \cR] \brho \;\; \qquad \quad
\quad & =  [\cD \cR] \br,
\end{align}
where $\bflambda \in \mathbb{R}^{2m}$ is a vector of Lagrange
multipliers.

In the numerical results presented in section \ref{sec:numres} we
consider smooth and piecewise constant resistivities, and choose ${W}$
in (\ref{eqn:h1w}) as follows.  For smooth resistivities we simply
take ${W}=I$, so that (\ref{eqn:h1w}) is a regular discrete $H^1$
seminorm. For discontinuous resistivities we would like to minimize
the total variation of the resistivity. This does not allow an
explicit computation of $\brho$, so we make a compromise and use the
weights introduced in \cite{abubakar2008}. The matrix ${W}$ is
diagonal with entries
\begin{equation}
{w}_{j} = \left(  ([\widetilde{D} \, \br]_j)^2 + \phi(\br)^2 \right)^{-1}, 
\quad j=1,\ldots,N-1,
\label{eqn:whatdiag}
\end{equation}
where $\phi(\br)$ is proportional to the misfit for the current iterate
\begin{equation}
\phi(\br) = C_\phi \| \cQ(d(\,\cdot\,)) - \cR(\br) \|_2,
\end{equation}
and $C_\phi$ is some constant, set to $1/(2m^2)$ in the numerical
examples in section \ref{sec:numres}.

To ensure that $A(\br)$ corresponds to a discretization of an elliptic
operator, we need positive entries in $\br$. This can be done with a
logarithmic change of coordinates, which transforms the optimization
problem to an unconstrained one. However, in our numerical experiments
we observed that if $m$ is sufficiently small so that for the given
data $d(t)$ all the entries of $\cQ(d(\,\cdot\,)) \in \mathbb{R}^{2m}$ 
are positive, then the Gauss-Newton updates of $\br$ remain positive as
well. Thus, the logarithmic change of coordinates for $\br$ in not
needed in our computations.

\subsection{The inversion algorithm for one dimensional media}
\label{subsec:invalg}

Here we present the summary of the inversion algorithm. The details 
of the computation of $\cR(\br)$ and its Jacobian $\cD \cR$ can be 
found in Appendix \ref{app:compjac}.

The inputs of the inversion algorithm
are the measured data $d(t)$ and a guess value of $m$. This $m$
reflects the expected resolution of the reconstruction, and may need
to be decreased depending on the noise level. To compute the estimate
$\br^\star$ of $\br^{\mbox{\tiny true}}$, perform the following steps:
\begin{itemize}
\item [1.]  Define the interpolation nodes $\widetilde{s}$ via
  (\ref{eqn:sigmatilde}).  Using the multipoint Pad\'{e} scheme from
  section \ref{subsec:ratinterp} compute $( \kappa^{\star}_j,
  \widehat{\kappa}_j^{\star} )_{j=1}^{m}$ using the data $d(t)$.
\item[2.] If for some $j$ either $\kappa^{\star}_j \leq 0$ or
  $\widehat{\kappa}^{\star}_j \leq 0$, decrease $m$ to $m-1$ and
  return to step 1.

  If all $( \kappa^{\star}_j, \widehat{\kappa}^{\star}_j )_{j=1}^{m}$
  are positive, fix $m$ and continue to step 3.
\item[3.] Define the vector of logarithms
 $\bl^\star = (\log \kappa_1^{\star}, \ldots, \log \kappa_m^{\star}, 
 \log \widehat{\kappa}_1^{\star}, \ldots, \log \widehat{\kappa}_m^{\star})^T$.
\item[4.] Choose an initial guess $\br^{ (1)} \in \mathcal{R}^N_+$ and
  the maximum number $n_{_{GN}}$  of Gauss-Newton iterations.
\item[5.] For $p=1,\ldots,n_{_{GN}}$ perform:
\item[] 
\begin{itemize}
\item[5.1.] For the current iterate $\br^{ (p)}$ compute the
  mapping \[\cR(\br^{ (p)}) = \{( \log \kappa^{ (p)}_j, \log
  \widehat{\kappa}_j^{ (p)} )\}_{j=1}^{m} \] and its
  Jacobian \[\cD \cR^{ (p)} = \cD \cR(\br^{ (p)})\] as
  explained in Appendix \ref{app:compjac}.
\item[5.2.] Define the vector of logarithms \[\bl^{ (p)} = (\log
  \kappa_1^{ (p)}, \ldots, \log \kappa_m^{ (p)}, \log
  \widehat{\kappa}_1^{ (p)}, \ldots, \log
  \widehat{\kappa}_m^{ (p)})^T.\]
\item[5.3.] Compute the step 
$$\bfrho^{(p)} = - \left( \cD \cR^{ (p)} \right)^\dagger 
(\bl^{ (p)} - \bl^\star).$$
\item[5.4.] Choose the step length $\alpha^{ (p)}$ and compute the
 Gauss-Newton update \[\br^{ GN} = \br^{ (p)} + \zeta^{ (p)}
  \bfrho^{(p)}.\]
\item[5.5.] Compute the weight $W$ using (\ref{eqn:whatdiag}) 
with $\br^{GN}$ or $W = I$.
\item[5.6.] Solve for the next iterate $\br^{ (p+1)}$ from
the linear system
\begin{equation}
\begin{bmatrix}
\widetilde{D}^T W \widetilde{D} & (\cD \cR^{ (p)})^T \\
\cD \cR^{ (p)} & 0
\end{bmatrix}
\begin{bmatrix}
\br^{ (p+1)} \\ \bflambda^{ (p)}
\end{bmatrix} =
\begin{bmatrix}
0 \\ (\cD \cR^{ (p)}) \br^{GN}
\end{bmatrix}
\label{eqn:reglinsolve}
\end{equation}
\end{itemize}
\item[6.] The estimate is 
\[
\br^\star = \br^{ (n_{_{GN}}+1)}.
\]
\end{itemize}

Let us remark that since the non-linear preconditioner is an
approximation of the identity in the sense explained in section
\ref{subsec:optim}, most of the nonlinearity of the problem is resolved
by the rational interpolation in step 1. Thus, we may start with a
poor initial guess in step 4 and still obtain good
reconstructions. Moreover, the number $n_{_{GN}}$ of Gauss-Newton
iterations may be kept small.  In the numerical results presented in
section \ref{sec:numres} we take the number of iterations
$n_{_{GN}} = 5$ for medium contrast resistivities and $n_{_{GN}} = 10$
for the high contrast case. In general, any of the standard stopping 
criteria could be used, such as the stagnation of the residual. 


In order to simplify our implementation we set the step length
$\alpha^{(p)} = 1$ in step 5.4.  However, choosing $\alpha^{(p)}$ 
adaptively with a line search procedure may be beneficial, 
especially for high contrast resistivities.

While the Jacobian $\cD \cR^{ (p)}$ is well-conditioned, the
system (\ref{eqn:reglinsolve}) may not be. To alleviate this problem,
instead of solving (\ref{eqn:reglinsolve}) directly, we may use a
truncated singular value decomposition to obtain a regularized
solution. Typically it is enough to discard just one component
corresponding to the smallest singular value as we do in the 
numerical experiments.

\section{Two dimensional inversion}
\label{sec:2D}

Unlike the one dimensional case, the inverse problem in two dimensions
is formally overdetermined. The unknown is the resistivity $r(x)$
defined on $\Omega \subset \mathbb{R}^2$, and the data are three
dimensional.  One dimension corresponds to time and the other two come
from the source (initial condition) and receiver locations on
$\mathcal{B}_A$. The model reduction inversion framework described in
section \ref{subsec:invalg} applies to a formally determined problem. We
extend it to two dimensions by constructing separately reduced models
for certain data subsets.  Each model defines a mapping $\cR_j$ that
is similar to $\cR$ in one dimension, where $j = 1, \ldots, N_d$ is
the index of the data set. The maps $\cR_j$ are coupled by their
dependence on the resistivity function $r(x)$, and are all taken into
account in inversion as explained below.

Let $u_o^{(j)}(x)$ be the initial condition for a source that is
compactly supported on a segment (interval) $\mathcal{J}_j$ of
$\mathcal{B}_A$. We model it for simplicity with the indicator
function $\mathds{1}_j$ of $\mathcal{J}_j$ and write
\begin{equation}
u_o^{(j)}(x) = \mathds{1}_j(x^\parallel)\delta(x^\perp).
\end{equation}
Here $x = (x^\parallel,x^\perp)$, with $x^\parallel$ the arclength on
$\mathcal{B}_A$ and $x^\perp$ the normal coordinate to the boundary,
which we suppose is smooth. The semidiscretized version of equation
(\ref{eq:3d}) on a grid with $N$ points is
\begin{equation}
\frac{\partial {\bf u}^{(j)}(t)}{\partial t} = A(\br) {\bf
  u}^{(j)}(t),
\label{eqn:semidfwd2D}
\end{equation}
with $\br \in \mathbb{R}^N$ the vector of discrete samples of the
resistivity.  Equation (\ref{eqn:semidfwd2D}) models a dynamical
system with response matrix $y_{kj}(t,\br)$ defined by the restriction
of the solution
\begin{equation}
{\bf u}^{(j)}(t) = e^{A(\br)t} {\bf u}_o^{(j)}
\end{equation}
to the support $\mathcal{J}_k$ of the $k-$th receiver. We take for
simplicity the same model of the sources and receivers, with support
on the disjoint boundary segments $\mathcal{J}_j$ of the accessible
boundary.  Thus, if we let ${\bf b}^{(j)} \in \mathbb{R}^N$ be the
measurement vector corresponding to the $j-$th source or receiver, we
can write the time domain response as
\begin{equation}
y_{kj}(t;\br) = {\bf b}^{(k)^T} e^{A(\br)t} {\bf b}^{(j)}.
\label{eqn:symy2D}
\end{equation}
The diagonal of this matrix is the high dimensional extension of
(\ref{eqn:symy}). The matrix valued transfer function is the Laplace
transform of (\ref{eqn:symy2D}), 
\begin{equation}
Y_{kj}(s; \br) = \int_{0}^{+\infty} y_{kj}(t; \br) e^{-st}dt =
{\bf b}^{(k)^T} (sI-A(\br))^{-1} {\bf b}^{(j)}, \quad
s>0.
\label{eqn:Gtransfer2D}
\end{equation}

In model order reduction we obtain a reduced model with rational
transfer function $Y_{kj,m}(s;\br)$ that approximates
(\ref{eqn:Gtransfer2D}). The reduced model is constructed separately
for each receiver-source pair $(k,j)$. It is defined by an $m \times
m$ symmetric and negative definite matrix $A_m^{(k,j)}(\br)$ with $m
\ll N$, and measurement vectors ${\bf b}_m^{(k)}$ and ${\bf
  b}_m^{(j)}$.  The transfer function
\begin{equation}
Y_{kj,m}(s; \br) = {\bf b}_m^{(k)^T} \left(s I_m -
A_m^{(k,j)}\right)^{-1} {\bf b}^{(j)} = \sum_{l=1}^m
\frac{c_l^{(k,j)}}{s + \theta_l^{(k,j)}},
\label{eqn:rat2D}
\end{equation}
has poles at the eigenvalues $-\theta_l^{(k,j)}$ of $A_m^{(k,j)}$, for
eigenvectots ${\bf z}_l^{(k,j)}$, and residues
\begin{equation}
c_l^{(k,j)} = \left({\bf b}_m^{(k)^T} {\bf
  z}_l^{(k,j)}\right)\left({\bf b}_m^{(j)^T} {\bf z}_l^{(k,j)}\right).
\label{eqn:rat2D1}
\end{equation}

We are interested in the continued fraction representation of
$Y_{kj,m}(s;\br)$, in particular its coefficients $\left\{\left(\widehat
\kappa_l^{(k,j)},\kappa_l^{(k,j)}\right)\right\}_{l=1}^m$ that define
the preconditioner mapping in our inversion approach. These
coefficients are guaranteed to be positive as long as the residues in
(\ref{eqn:rat2D1}) satisfy $c_l^{(k,j)} \ge 0$. This is guaranteed to
hold for the diagonal of (\ref{eqn:rat2D}), because
\begin{equation}
c_l^{(j,j)} = \left({\bf b}_m^{(j)^T} {\bf z}_l^{(j,j)}\right)^2.
\label{eqn:rat2D2}
\end{equation}
Thus, we construct the reduced models for the diagonal of the matrix
valued transfer function, and compute the model parameters
$\left\{\left(\widehat
\kappa_l^{(j,j)},\kappa_l^{(j,j)}\right)\right\}_{l=1}^m$ as in one
dimension. Such measurement setting has analogues in other types 
of inverse problems. For example, a similar setting in wave inversion
is the backscattering problem, where the scattered wave field is 
measured in the same direction as the incoming wave.

If we have $N_d$ boundary segments $\mathcal{J}_j$ at $\mathcal{B}_A$,
we define
\begin{equation}
\cR_j(\br) = \cQ\left(y_{jj}(\cdot; \br)\right) = \left(\log \widehat
\kappa_l^{(j,j)},\log \kappa_l^{(j,j)}\right)_{l=1}^m, \qquad j = 1,
\ldots, N_d,
\end{equation}
as in one dimension, using the chain of mappings
(\ref{eqn:kappachain}). The resistivity is estimated by the solution
of the optimization problem
\begin{equation}
\br^\star = \mathop{\mbox{arg min}}\limits_{\br \in \mathbb{R}^N_+}
\frac{1}{2}\sum_{j=1}^{N_d} \left\| \cQ(d_j(\,\cdot\,)) - \cR_j(\br)) \right\|_2^2,
\label{eqn:roptim2D}
\end{equation}
where $d_j(t)$ is the data measured at the receiver supported on
$\mathcal{J}_j$, for the $j-th$ source excitation. Note that the
sum in (\ref{eqn:roptim2D}) couples all the sources/receivers 
together in a single objective functional. 

Another question that we need to address is the choice of matching
conditions. For simplicity we match the moments of $Y_{jj}(s)$ at a 
single interpolation node $\widetilde{s}$. This yields the matching 
conditions
\begin{equation}
\left. \frac{\partial^k Y_{jj,m}}{\partial s^k} \right|_{s = \widetilde{s}} = 
\left. \frac{\partial^k Y_{jj}}{\partial s^k} \right|_{s = \widetilde{s}},
\quad k = 0,1,\ldots,2m-1,
\label{eqn:2Dmatchcond}
\end{equation} 
and rational Krylov subspaces
\begin{equation}
\mathcal{K}_{m}^{(j)} (\widetilde{s}) = \mbox{span} \left\{ 
(\widetilde{s} I - A)^{-1} \bb^{(j)}, \ldots,
(\widetilde{s} I - A)^{-m} \bb^{(j)}\right\}, \quad
j = 1,\ldots,N_d.
\end{equation} 
Then the only parameter to determine is the node $\widetilde{s} > 0$.

Similarly to the one dimensional case we use the condition number of
the Jacobian $\cD \cR$ to determine the optimal choice of 
$\widetilde{s}$. Here the total Jacobian 
$\cD \cR \in \mathbb{R}^{(2 m N_d) \times N}$ is a matrix of 
$N_d$ individual Jacobians $\cD \cR_j \in \mathbb{R}^{2 m \times N}$
stacked together. It appears that for a fixed $j$ the sensitivity
functions for $\log \kappa_l^{(j,j)}$ and $\log \widehat{\kappa}_l^{(j,j)}$
closely resemble the spherical waves propagating from the
middle of $\mathcal{J}_j$. As the index $l$ increases, the wave 
propagates further away from $\mathcal{J}_j$. This behavior is 
illustrated quantitatively in section \ref{subsec:numset2} for an 
optimal choice of $\widetilde{s}$. The speed of propagation of
these sensitivity waves decreases as $\widetilde{s}$ grows. 
Obviously, to get the resolution throughout the whole domain we
would like all of $\Omega$ to be covered by the sensitivity waves,
which means higher propagation speed (smaller $\widetilde{s}$) is
needed. On the other hand, if the sensitivity waves propagate too 
far they reflect from the boundary $\cB$ which leads to poor 
conditioning of the Jacobian. The balance between these two 
requirements leads to an optimal choice of $\widetilde{s}$,
which we determine experimentally in the numerical example in 
section \ref{subsec:numset2}.

\section{Numerical results}
\label{sec:numres}

We assess the performance of the inversion algorithm with numerical
experiments.  To avoid committing the inverse crime we use different
grids for generating the data and for the solution of the inverse
problem.  We describe the setup of the numerical simulations and
present the inversion results in section \ref{subsec:numset} for one
dimensional media, and in section \ref{subsec:numset2} for two
dimensional media.

\subsection{Numerical experiments in one dimension}
\label{subsec:numset}
We use a fine grid with $N_f = 299$ uniformly spaced nodes to simulate
the data $d(t)$, and a coarser grid with $N = 199$ uniformly spaced
nodes in the inversion.  

The first term $y(t; \br^{\mbox{\tiny true}})$ in
(\ref{eqn:semiddata}) is approximated by solving the semi-discrete
forward problem (\ref{eqn:semidfwd}) with an explicit forward Euler
time stepping on a finite time interval $[0, T]$, where $T = 100$. The
time step is $h_T = 10^{-5}$.  We denote by $\by$ the vector of length
$N_T = T/h_T$, with entries given by the numerical approximation of
$y(t;\br^{\mbox{\tiny true}})$ at the time samples $t_j = j
h_T$. Since even in the absence of noise there is a systematic error
$\mathcal{N}^{ (s)}(t)$ coming from the numerical approximation
of the solution of (\ref{eqn:semidfwd}), we write
\[
\by = (y_1, \ldots, y_{_{N_T}})^T, \quad y_j = y(t_j;\br^{\mbox{\tiny
    true}}) + \mathcal{N}^{(s)}(t_j), \quad j = 1, \ldots, N_T.
\]
We also define the vector $\mathcal{N}^{ (n)} \in
\mathbb{R}^{N_T}$ that simulates measurement noise using the
multiplicative model
\begin{equation}
\boldsymbol{\mathcal{N}}^{ (n)} = \epsilon \,
\mbox{diag}(\chi_1,\ldots,\chi_{N_T}) \by,
\end{equation}
where $\epsilon$ is the noise level, and $\chi_k$ are independent
random variables distributed normally with zero mean and unit standard
deviation.  The data vector ${\bf d} \in \mathbb{R}^{N_T}$ is 
\begin{equation}
{\bf d} = {\bf y} + \boldsymbol{\mathcal{N}}^{ (n)},
\end{equation}
and we denote its components by $d_j$, for $j = 1, \ldots, N_T$.
Such noise model allows for a simple estimate for a 
signal-to-noise ratio
\begin{equation}
\frac{\|\bd\|_2}{\| \boldsymbol{\mathcal{N}}^{ (n)} \|_2} \approx \frac{1}{\epsilon}.
\end{equation}

The inversion algorithm described in section \ref{subsec:invalg}
determines at step 2 the size $m$ of the reduced model for different
levels of noise.  The larger $\epsilon$, the smaller $m$. The values
of $\epsilon$ and $m$ used to obtain the results presented here are
given in Table \ref{tab:noiselevel}.

\begin{table}[ht!]
\begin{center}
\caption{Reduced model sizes $m$ used for various noise levels $\epsilon$.}
\begin{tabular}{c|c|c|c|c}
$\epsilon$ & $5\cdot10^{-2}$ & $5\cdot10^{-3}$ & $10^{-4}$ & $0$ (noiseless) \\
\hline
$m$ & 3 & 4 & 5 & 6
\end{tabular}
\label{tab:noiselevel}
\end{center}
\end{table}

The transfer function and its derivative at the interpolation points
are approximated by taking the discrete Laplace transform of the
simulated data
\begin{eqnarray}
Y(\widetilde{s}_j) & \approx & h_T \sum_{k=1}^{N_T} d_k
e^{- \widetilde{s}_j t_k},
\label{eqn:Gdlaplace}\\
Y'(\widetilde{s}_j) & \approx & - h_T \sum_{k=1}^{N_T}
d_k t_k e^{- \widetilde{s}_j t_k}.
\label{eqn:Gdplaplace}
\end{eqnarray}

\begin{figure}[ht!]
\begin{center}
\includegraphics[width=0.32\textwidth]{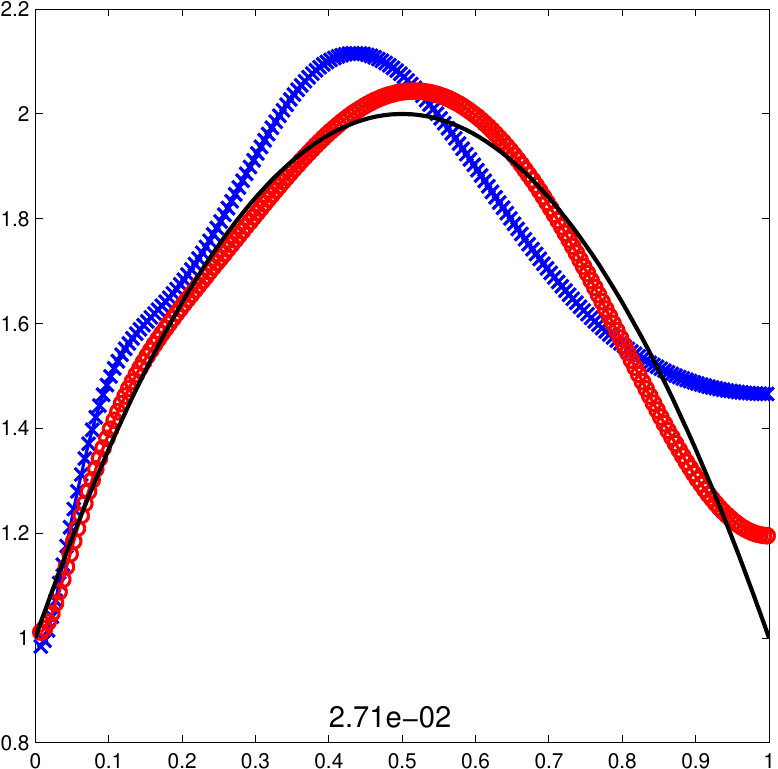}
\includegraphics[width=0.32\textwidth]{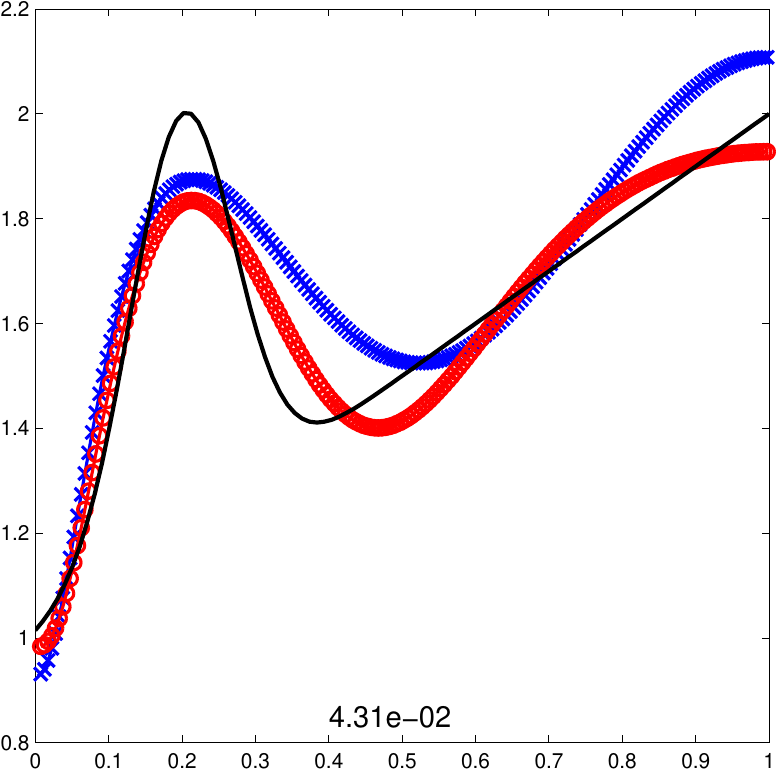}
\includegraphics[width=0.32\textwidth]{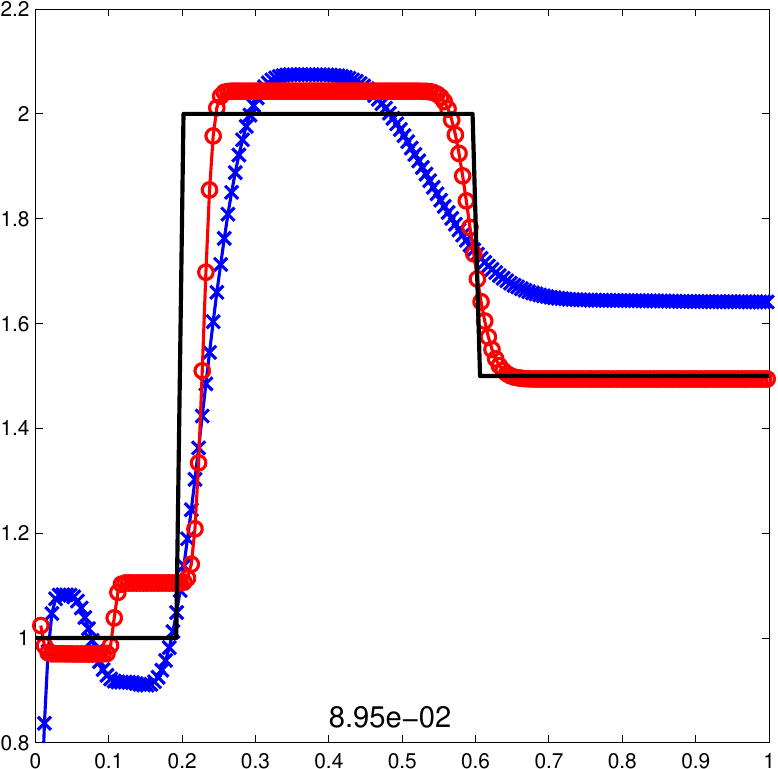}
\includegraphics[width=0.32\textwidth]{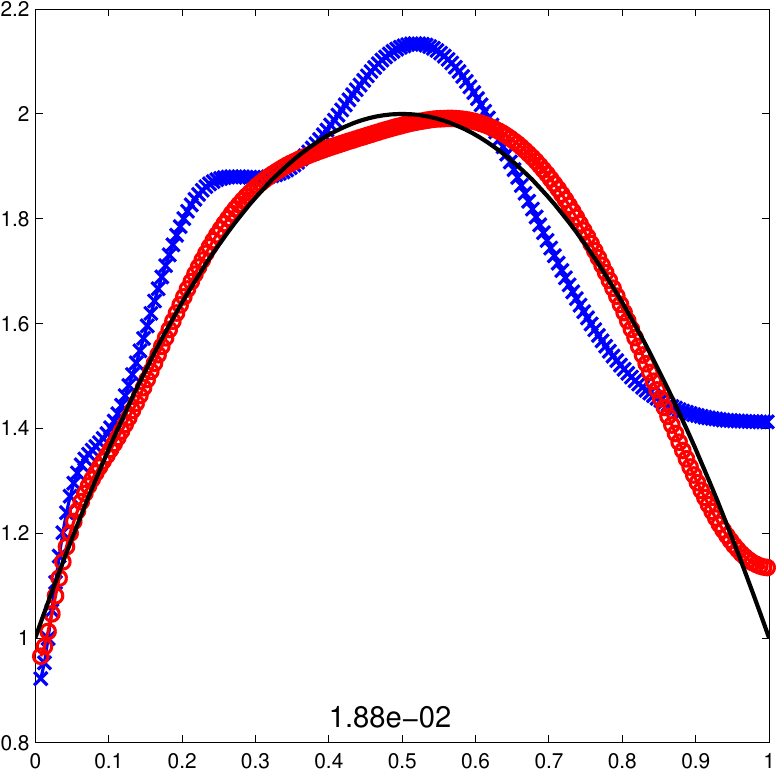}
\includegraphics[width=0.32\textwidth]{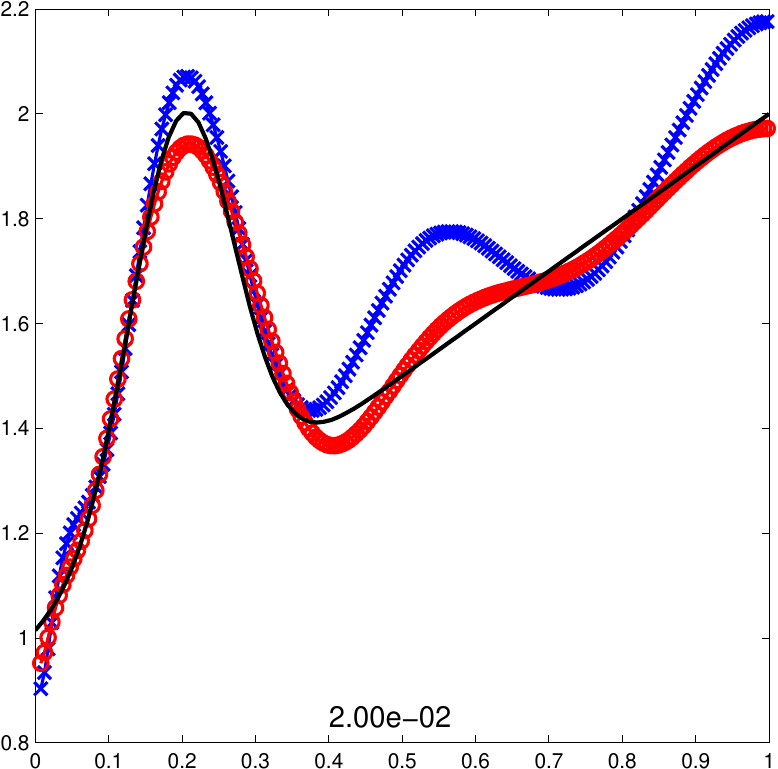}
\includegraphics[width=0.32\textwidth]{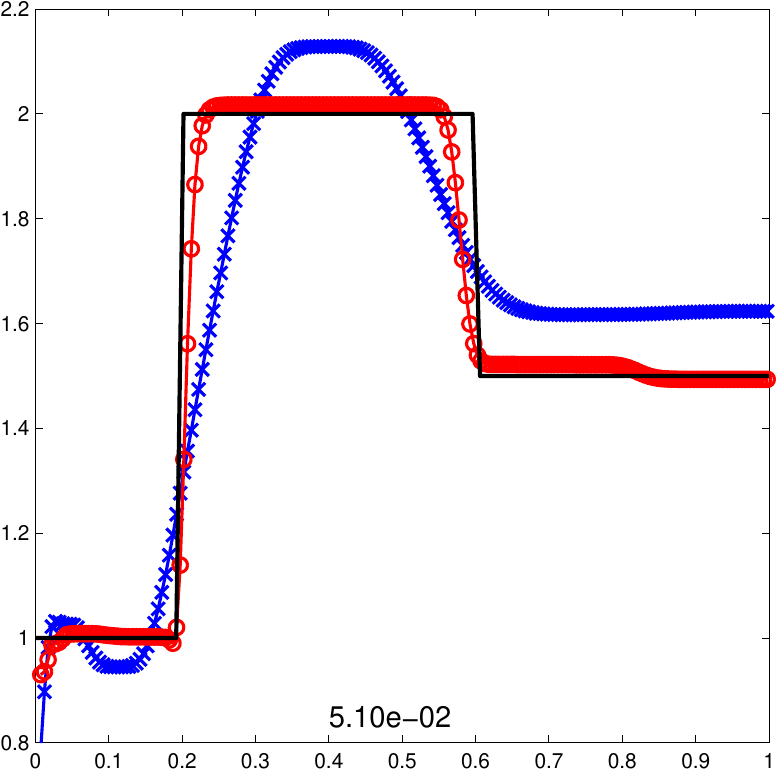}
\includegraphics[width=0.32\textwidth]{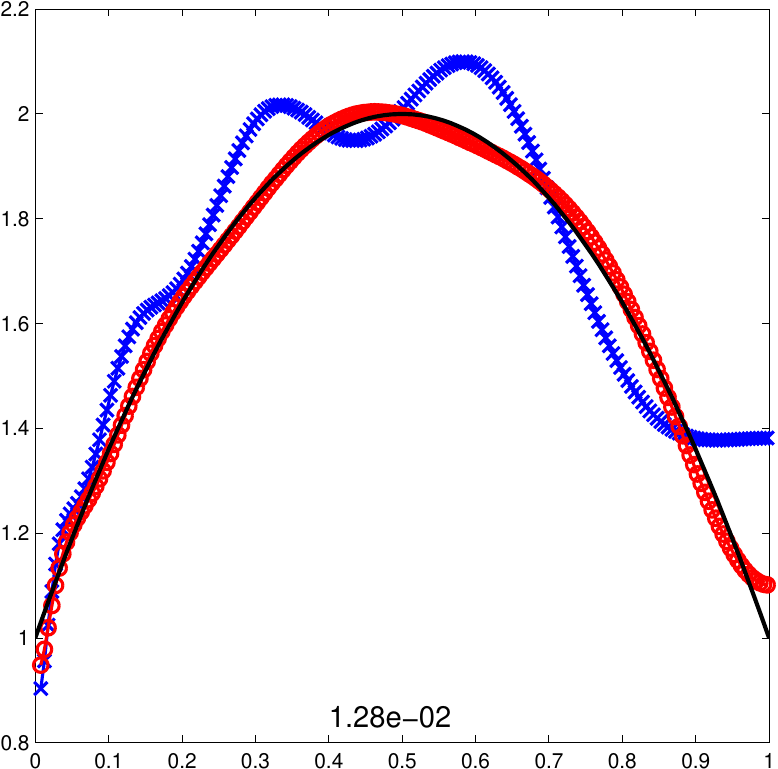}
\includegraphics[width=0.32\textwidth]{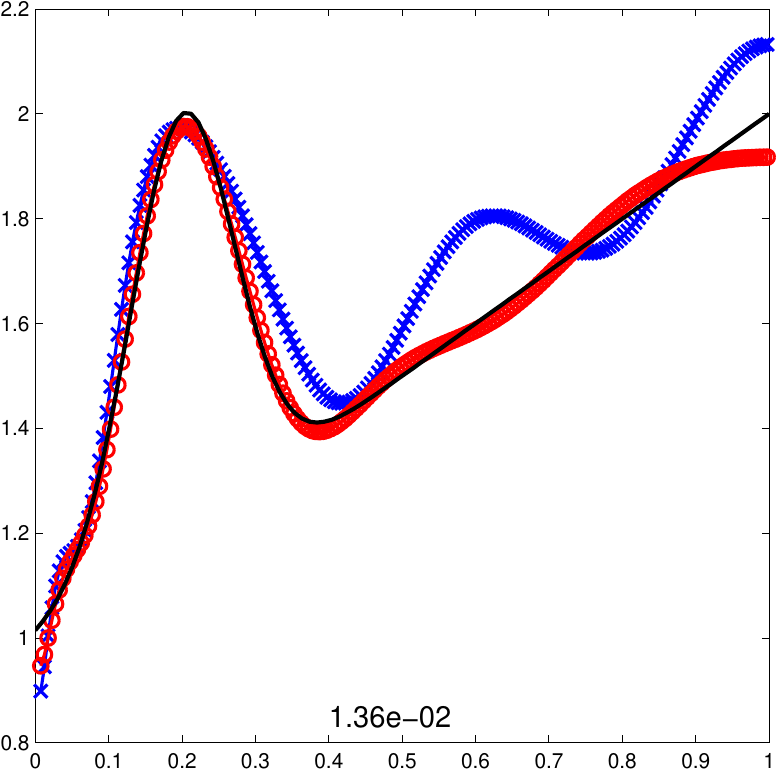}
\includegraphics[width=0.32\textwidth]{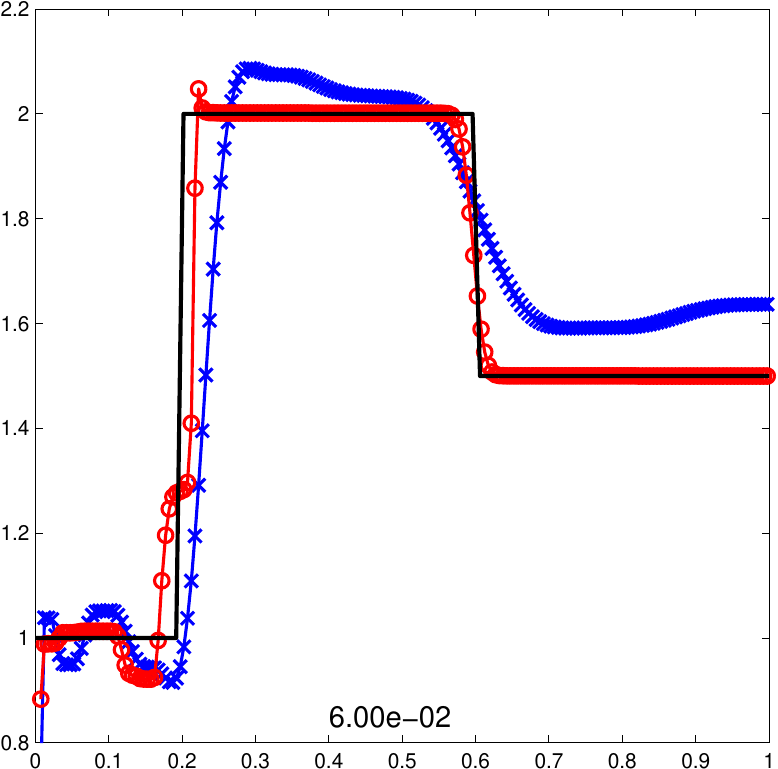}
\includegraphics[width=0.32\textwidth]{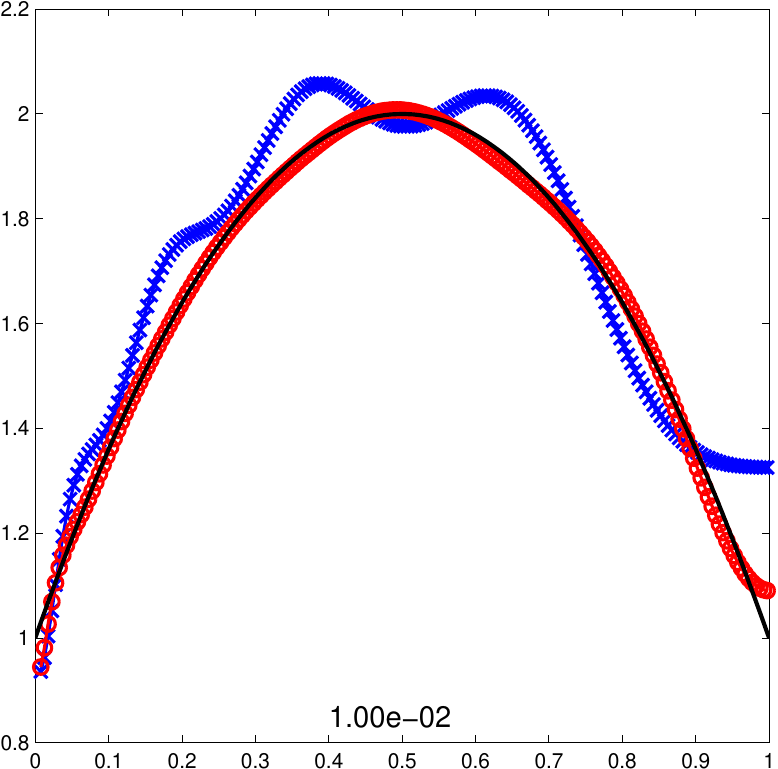}
\includegraphics[width=0.32\textwidth]{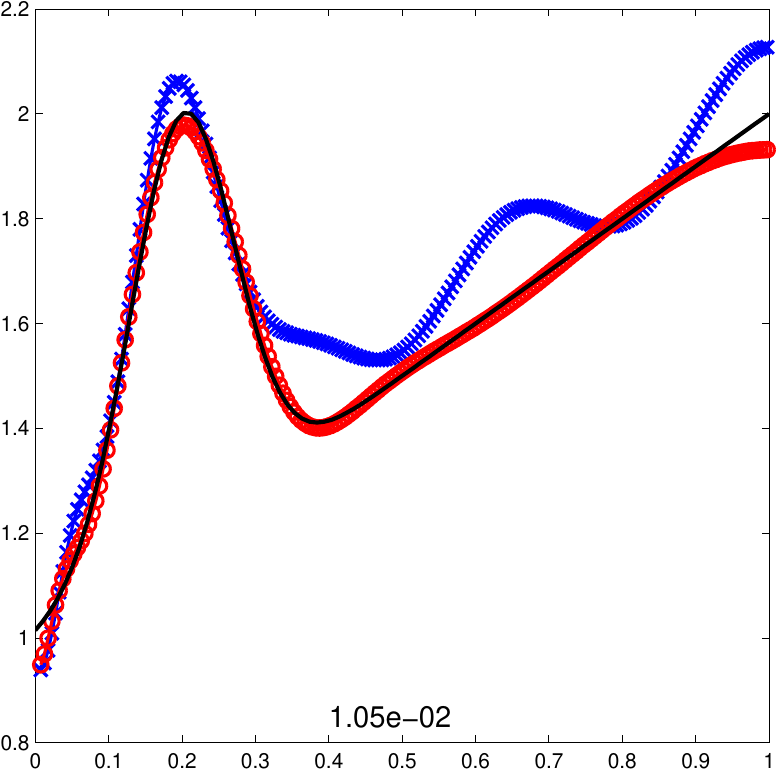}
\includegraphics[width=0.32\textwidth]{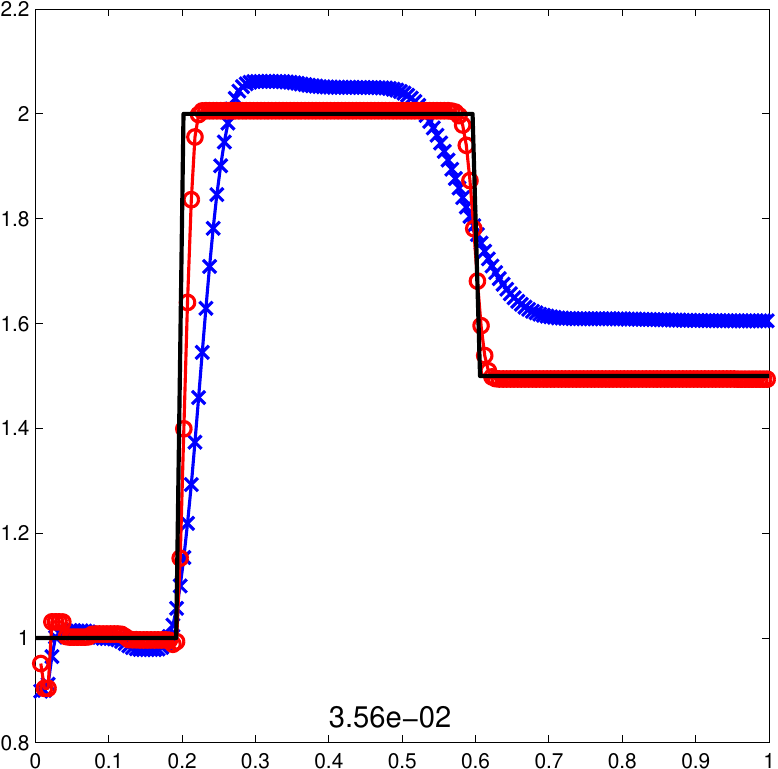}
\caption{Reconstructions of $r(x)$ (black solid line) after one 
(blue $\times$) and five (red $\circ$) iterations. True coefficient 
by column: left $r_Q$, middle $r_L$, right $r_J$.
Reduced model size from top row to bottom row $m=3,4,5,6$. 
Noise levels are from Table \ref{tab:noiselevel}.
The relative error $\mathcal{E}$ is printed at the bottom of the plots.}
\label{fig:numgn}
\end{center}
\end{figure}

To quantify the error of the reconstructions $\br^\star$ we use the
ratio of discrete $\ell_2$ norms
\begin{equation}
\label{eq:error}
\mathcal{E} = \frac{\| \br^\star - \br^{\mbox{\tiny true}} \|_2}
{\| \br^{\mbox{\tiny true}} \|_2}.
\end{equation}
While this measure of error is most appropriate for smooth
resistivities, it may be overly strict for the reconstructions in the
discontinuous case due to the large contribution of discontinuities.
However, even under such an unfavorable measure the inversion
procedure demonstrates good performance.

We show first the estimates of three resistivity functions of contrast
two. We consider two smooth resistivities
\begin{eqnarray}
&&r^{\mbox{\tiny true}}(x) = r_Q(x)  :=  2 - 4 \left( x - \frac{1}{2} \right)^2, \\
&&r^{\mbox{\tiny true}}(x) = r_L(x)  :=  0.8 e^{-100(x - 0.2)^2} + x + 1,
\label{eqn:rQL}
\end{eqnarray}
and the piecewise constant 
\begin{equation}
r^{\mbox{\tiny true}}(x) = r_J(x) := \left\{ 
\begin{tabular}{ll} 
1, & for $x<0.2$ \\
2, & for $0.2 \leq x \leq 0.6$ \\
1.5, & for $x > 0.6$ \\
 \end{tabular} \right.
\label{eqn:rJ}
\end{equation}
The results are displayed in Figure \ref{fig:numgn}, for various
reduced model sizes and levels of noise, as listed in Table
\ref{tab:noiselevel}. Each reconstruction uses its own realization of
noise. We use five Gauss-Newton iterations $n_{_{GN}} = 5$, and we
display the solution both after one iteration and after all five. The
initial guess is $r^{(1)}(x) \equiv 1$. It is far from the true
resistivity $r^{\mbox{\tiny true}}(x)$, and yet the inversion
procedure converges quickly. The features of
$r^{\mbox{\tiny true}}(x)$ are captured well after the first
iteration, but there are some spurious oscillations, corresponding to
the peaks of the sensitivity functions. A few more iterations of the 
inversion algorithm remove these oscillations and improve the
quality of the estimate of the resistivity. The relative error
(\ref{eq:error}) is indicated in each plot in Figure
\ref{fig:numgn}. It is small, of a few percent in all cases.

\begin{figure}[ht!]
\begin{center}
\begin{tabular}{ccc}
$m=4$, $\epsilon = 5\cdot10^{-3}$ & $m=5$, $\epsilon = 10^{-4}$ & $m=6$, $\epsilon = 0$ \\
\hskip-0.01\textwidth \includegraphics[width=0.32\textwidth]{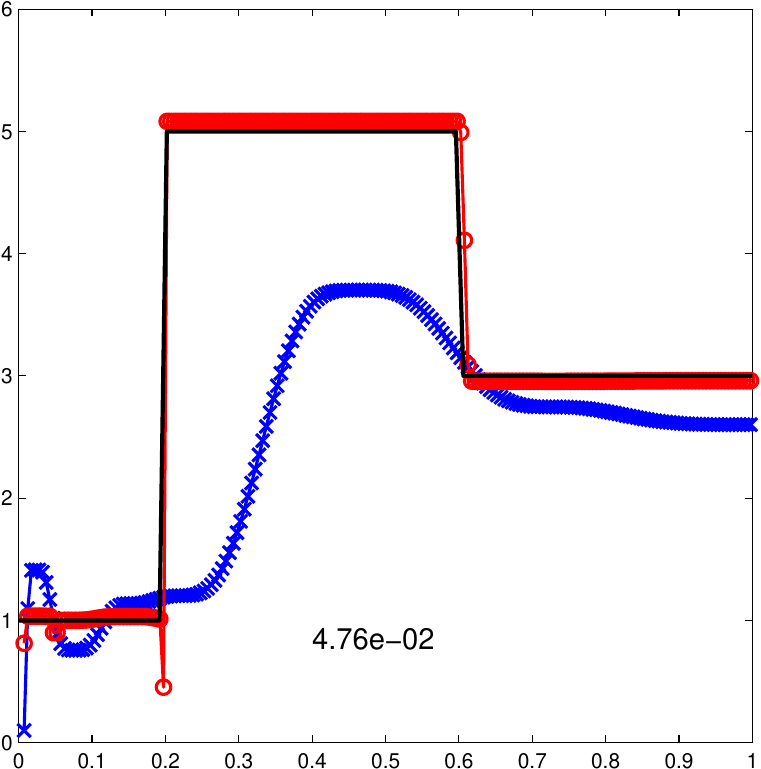} & 
\hskip-0.01\textwidth \includegraphics[width=0.32\textwidth]{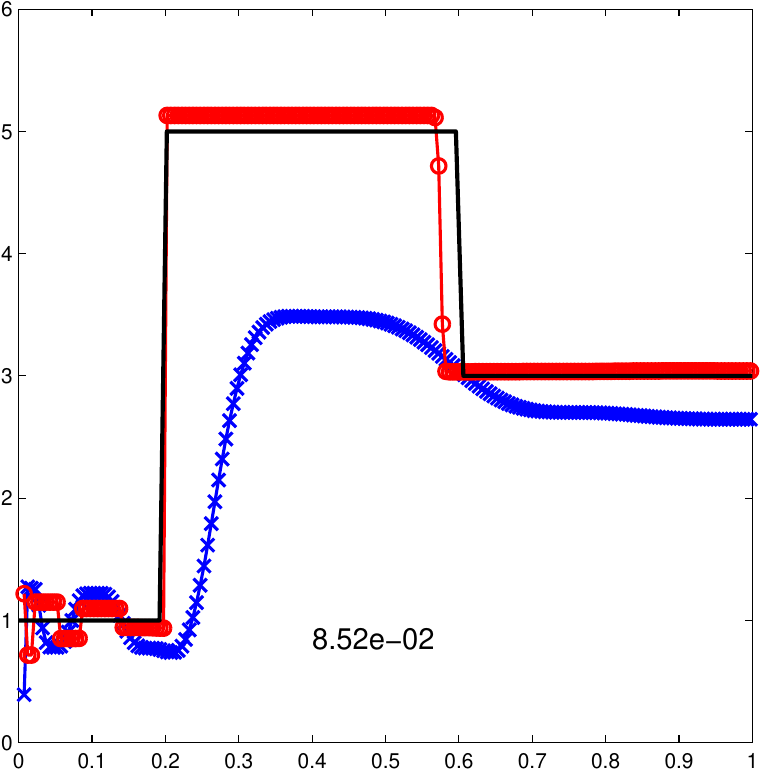} &
\hskip-0.01\textwidth \includegraphics[width=0.32\textwidth]{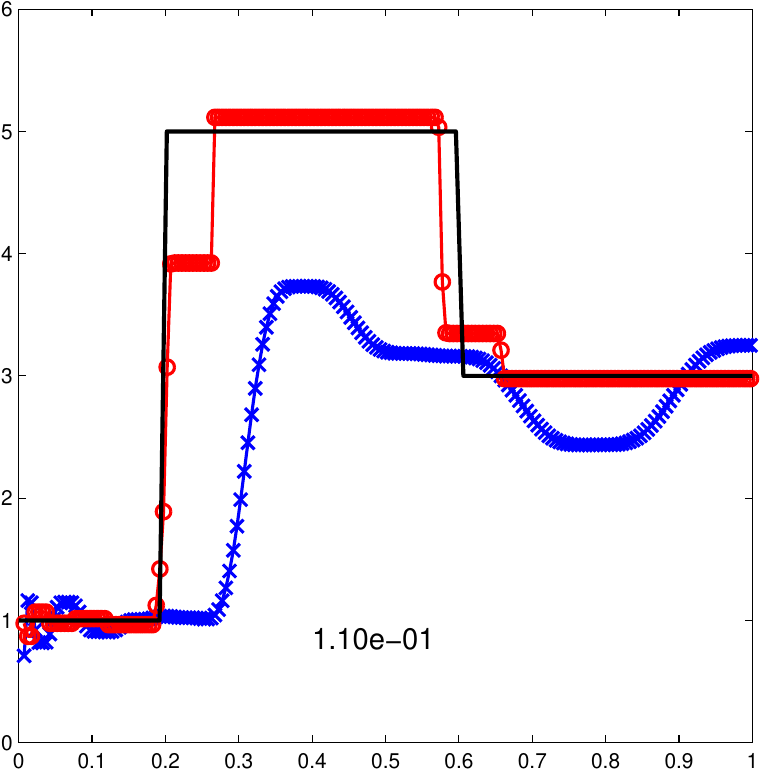} \\
\hskip-0.01\textwidth \includegraphics[width=0.32\textwidth]{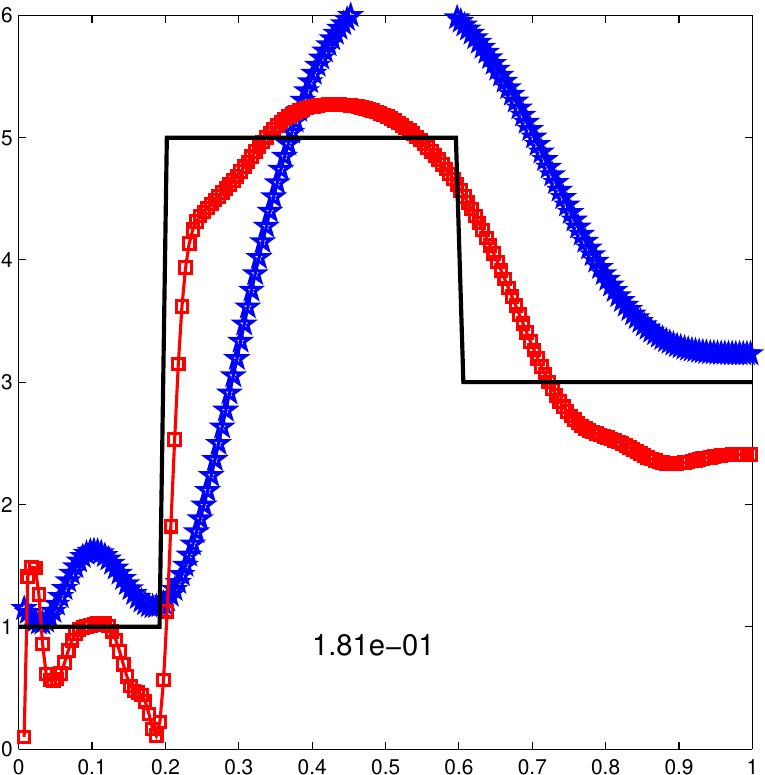} &
\hskip-0.01\textwidth \includegraphics[width=0.32\textwidth]{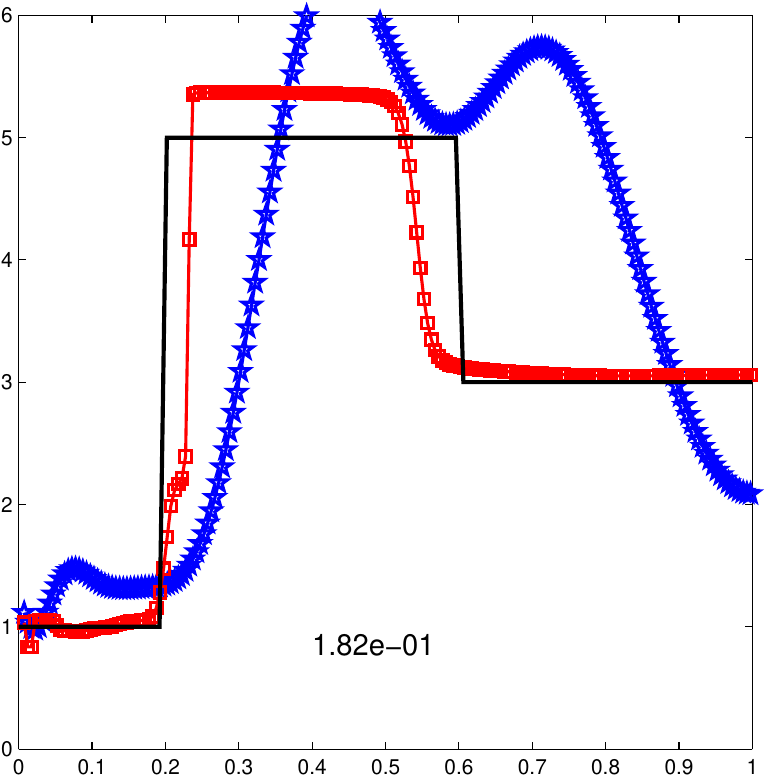} &
\hskip-0.01\textwidth \includegraphics[width=0.32\textwidth]{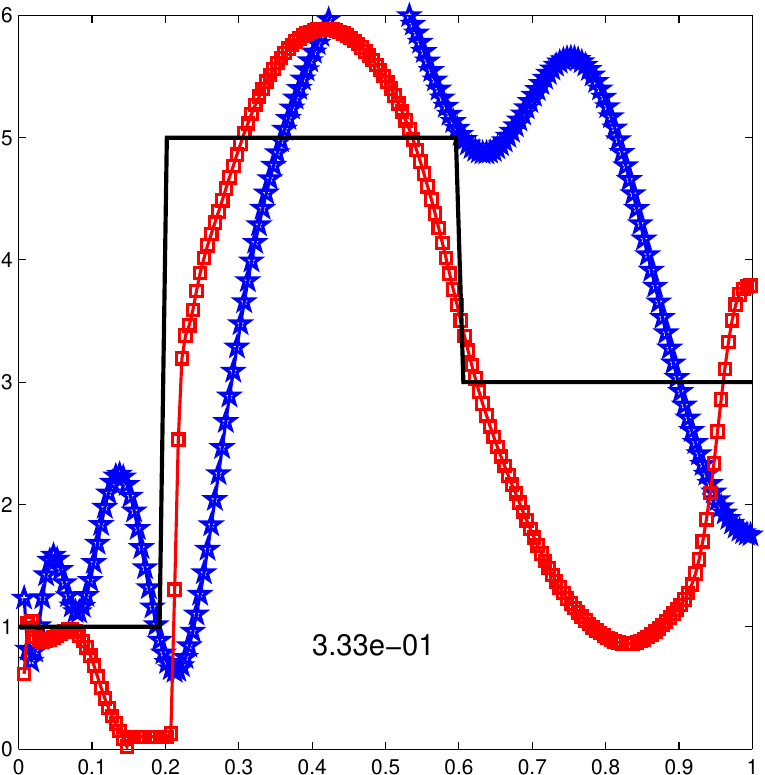} \\
\hskip-0.02\textwidth \includegraphics[width=0.33\textwidth]{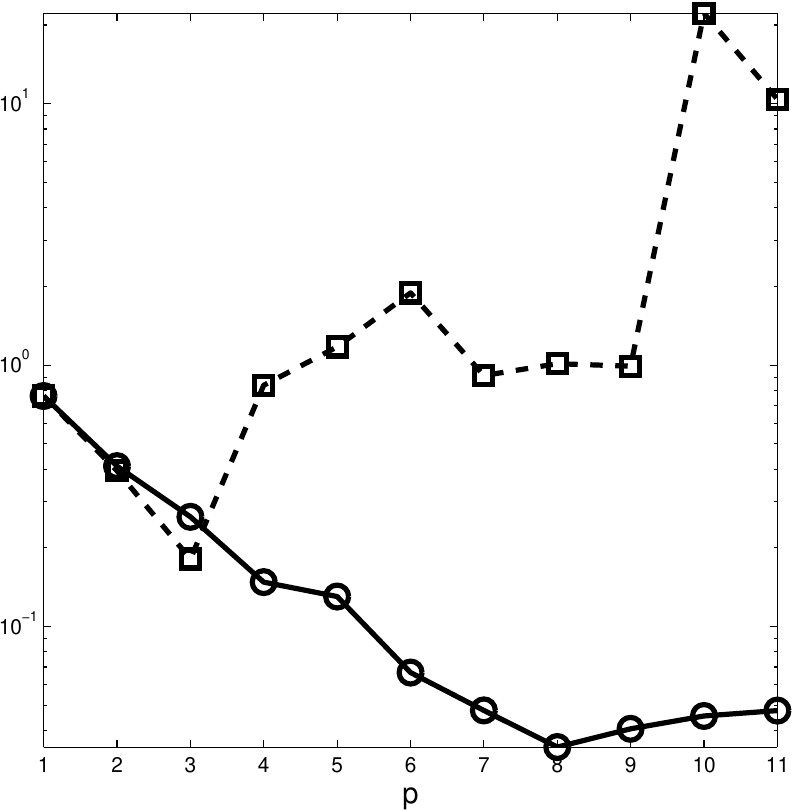} &
\hskip-0.02\textwidth \includegraphics[width=0.33\textwidth]{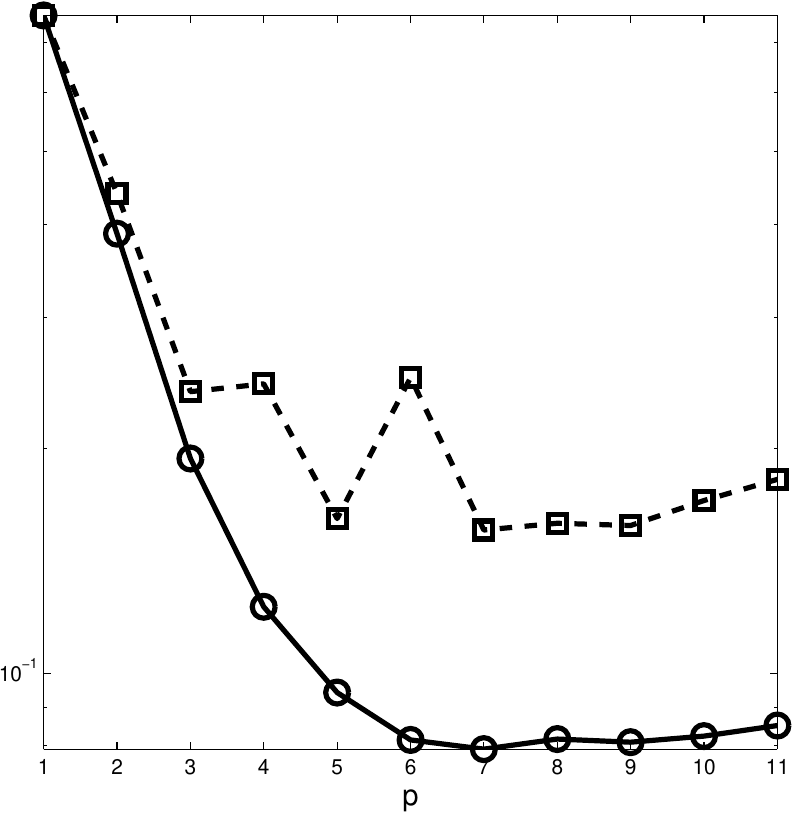} &
\hskip-0.02\textwidth \includegraphics[width=0.33\textwidth]{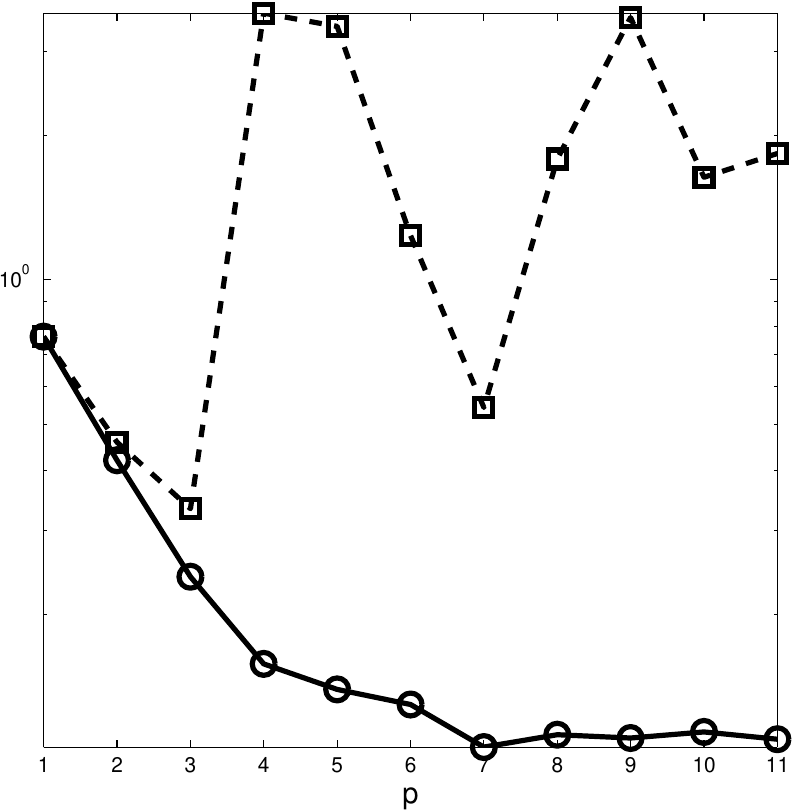}
\end{tabular}
\caption{Comparison between two preconditioners for high 
contrast piecewise constant resistivity $r_H(x)$ (black solid line)
after one (blue $\times$ and $\star$) and $n_{GN}$ (red $\circ$ and $\square$) 
iterations. 
Top row: reconstructions using $\cR$ ($n_{GN}=10$).
Middle row: reconstructions using $\cR_{_{c\theta}}$
($n_{GN}=2$ for $m=4,6$ and $n_{GN}=10$ for $m=5$). 
Bottom row: relative error $\cE$ versus the iteration number $p$
for reconstructions with $\cR$ (solid line with $\circ$) and 
$\cR_{_{c\theta}}$ (dashed line with $\square$).
The relative error $\cE$ is printed at the bottom of the reconstruction plots.
}
\label{fig:numhigh}
\end{center}
\end{figure}

We also observe in Figure \ref{fig:numgn} that the inversion method
regularized with the non-linear weight (\ref{eqn:whatdiag}) performs
well for the piecewise constant resistivity $r_J$. Without the 
regularization, the estimates have Gibbs-like oscillations near the 
discontinuities of $r^{\mbox{\tiny true}}(x)$. These oscillations are 
suppressed by the weighted discrete $H^1$ regularization.

In Figure \ref{fig:numhigh} we are comparing our inversion procedure
to an inversion approach like in \cite{dsz2012} that fits the poles
and residues $(\theta_j,c_j)_{j=1}^{m}$ instead of the continued
fraction coefficients. 
The comparison is done for the case of piecewise constant resistivity of
higher contrast
\begin{equation}
r^{\mbox{\tiny true}}(x) = r_H(x) := \left\{ 
\begin{tabular}{ll} 
1, & for $x<0.2$ \\
5, & for $0.2 \leq x \leq 0.6$ \\
3, & for $x > 0.6$ \\
\end{tabular} \right.
\end{equation}
A higher contrast case is chosen since for the low contrast the 
difference in performance between $\cR_{_{c\theta}}$ and $\cR$ is
less pronounced.

As expected, the reconstructions are better when we use the mapping
$\cR$. In fact, the algorithm based on $\cR_{_{c\theta}}$ diverges for
$m=4$ and $m=6$. Thus, for $m=4,6$ we plot the first and second
iterates only. The inversion based on the mapping $\cR$ converges in
all three cases and maintains the relative error well below $10\%$ for
$m=4,5$ and around $11\%$ for $m=6$.  The reconstruction plots in
Figure \ref{fig:numhigh} are complimented with the plots of the
relative error $\cE$ versus the Gauss-Newton iteration number
$p$. Note that even when the iteration with $\cR_{_{c\theta}}$ converges
($m=5$) the reconstruction with $\cR$ has smaller error ($8.5\%$
versus $18.2\%$).

Aside from providing a solution of higher quality our method is also
more computationally efficient since it does not require the solution
of the forward problem (\ref{eqn:semidfwd}) in time. While
constructing the orthonormal basis for the Krylov subspace
(\ref{eqn:ratkrylov}) requires a few linear solves with the shifted
matrix $A$, the number $m$ of such solves is small and thus it is
cheaper than the time stepping for (\ref{eqn:semidfwd}). For example,
the explicit time stepping to generate the data $\bd$ for the
numerical experiments above takes $275$ seconds, whereas
all the five Gauss-Newton iterations of our inversion algorithm takes less
than a second on the same machine.

\subsection{Numerical experiments in two dimensions}
\label{subsec:numset2}

\begin{figure}[t]
\begin{center}
\begin{tabular}{ccc}
& $\dfrac{\partial \log \kappa_l^{(j,j)}}{\partial r}_{\;_{\;}}$ & 
  $\dfrac{\partial \log \widehat{\kappa}_l^{(j,j)}}{\partial r}_{\;_{\;}}$ \\
$(l=1)$ &
\hskip-0.02\textwidth
\raisebox{-0.5\totalheight}{\includegraphics[width=0.44\textwidth]{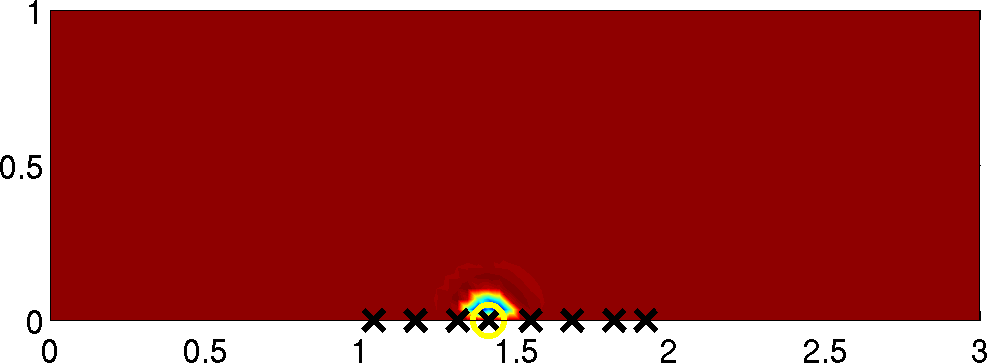}} &
\hskip-0.02\textwidth
\raisebox{-0.5\totalheight}{\includegraphics[width=0.44\textwidth]{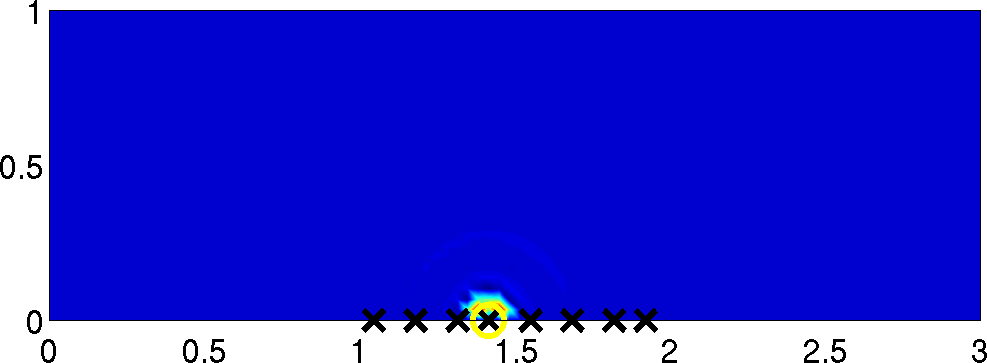}} \\
$(l=2)$ &
\hskip-0.02\textwidth
\raisebox{-0.5\totalheight}{\includegraphics[width=0.44\textwidth]{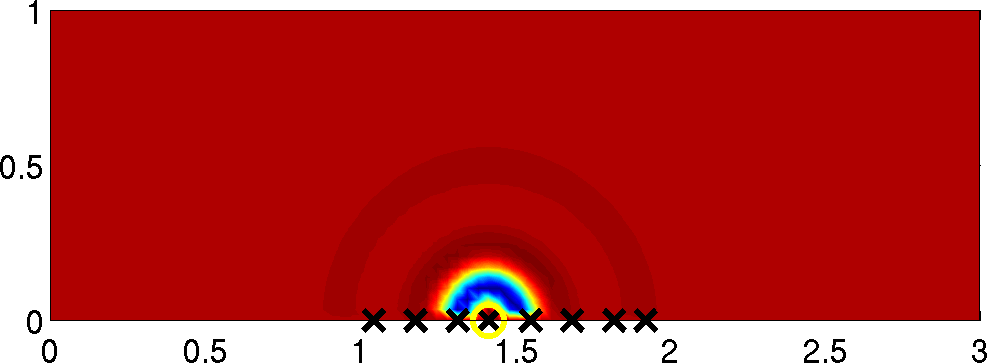}} & 
\hskip-0.02\textwidth
\raisebox{-0.5\totalheight}{\includegraphics[width=0.44\textwidth]{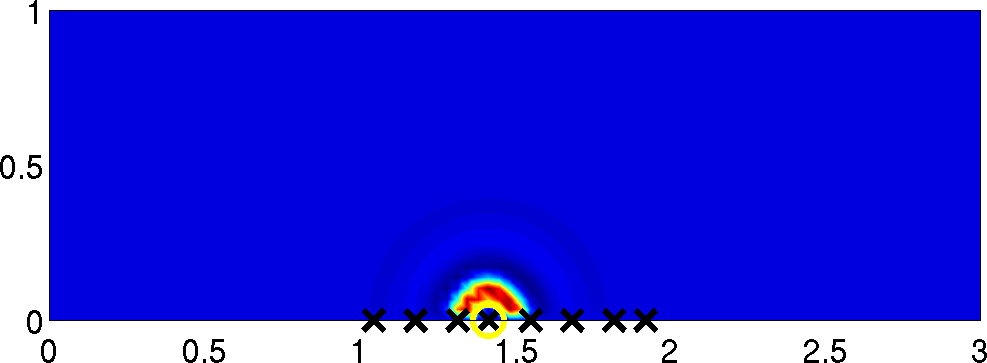}} \\
$(l=3)$ &
\hskip-0.02\textwidth
\raisebox{-0.5\totalheight}{\includegraphics[width=0.44\textwidth]{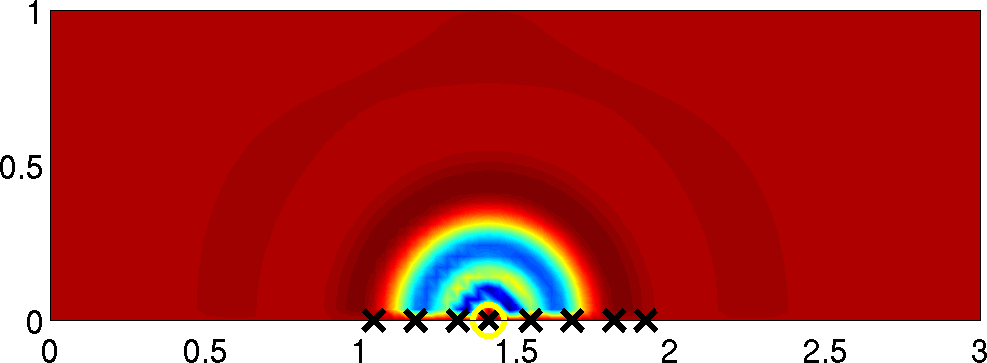}} & 
\hskip-0.02\textwidth
\raisebox{-0.5\totalheight}{\includegraphics[width=0.44\textwidth]{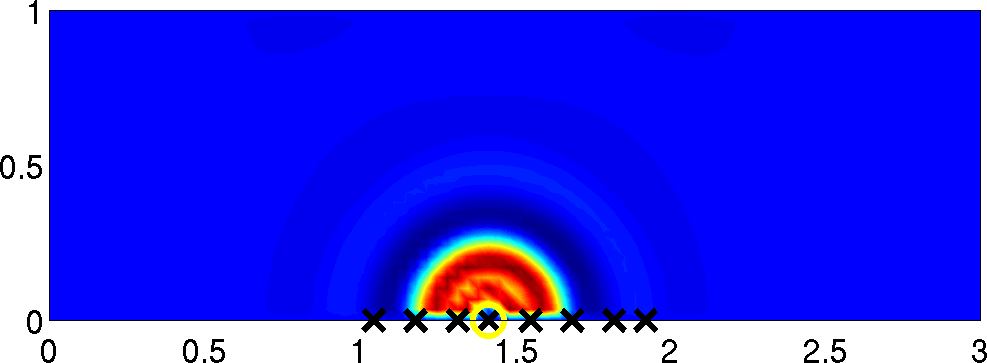}} \\
$(l=4)$ &
\hskip-0.02\textwidth
\raisebox{-0.5\totalheight}{\includegraphics[width=0.44\textwidth]{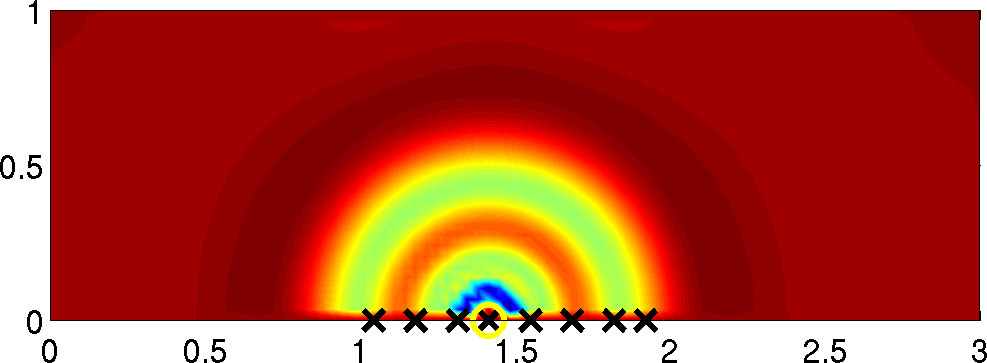}} & 
\hskip-0.02\textwidth
\raisebox{-0.5\totalheight}{\includegraphics[width=0.44\textwidth]{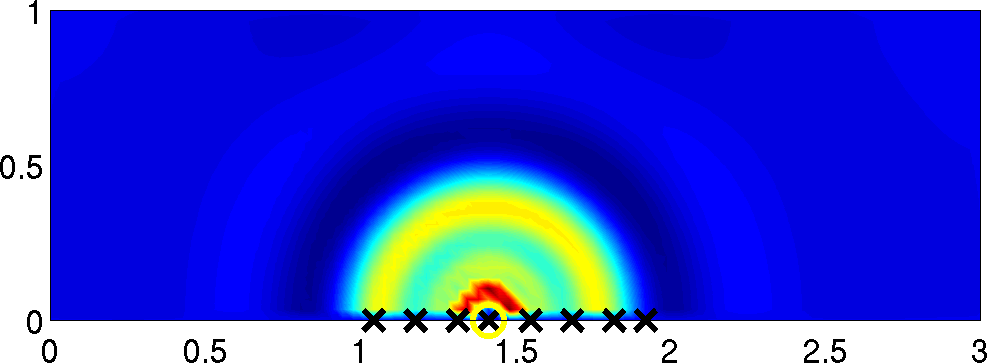}} \\
$(l=5)$ &
\hskip-0.02\textwidth
\raisebox{-0.5\totalheight}{\includegraphics[width=0.44\textwidth]{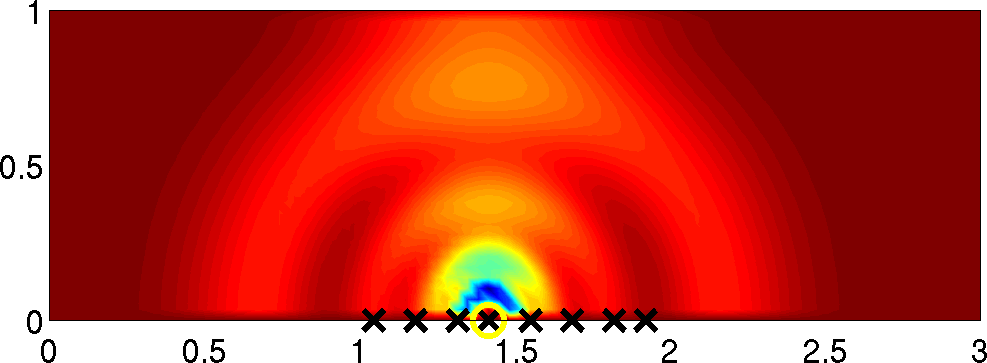}} & 
\hskip-0.02\textwidth
\raisebox{-0.5\totalheight}{\includegraphics[width=0.44\textwidth]{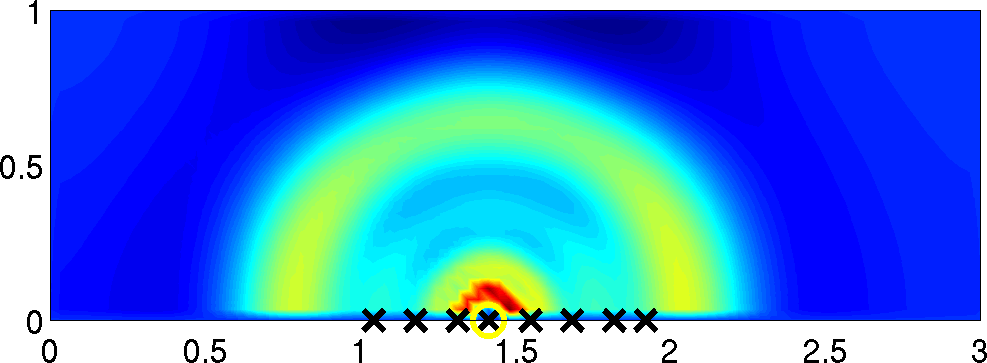}}
\end{tabular}
\caption{Sensitivity functions (rows of the Jacobian $\cD \cR_j$ indexed by $l$) for
the two dimensional uniform medium $r(x) \equiv 1$ in the rectangular 
domain $\Omega = [0,3] \times [0,1]$. Source/receiver index is $j=4$
out of a total of $N_d = 8$ (mid-points of sources/receivers $\cJ_j$ marked 
as black $\times$), $m = 5$.
}
\label{fig:2Dsens}
\end{center}
\end{figure}

\begin{figure}[t]
\begin{center}
\includegraphics[width=0.49\textwidth]{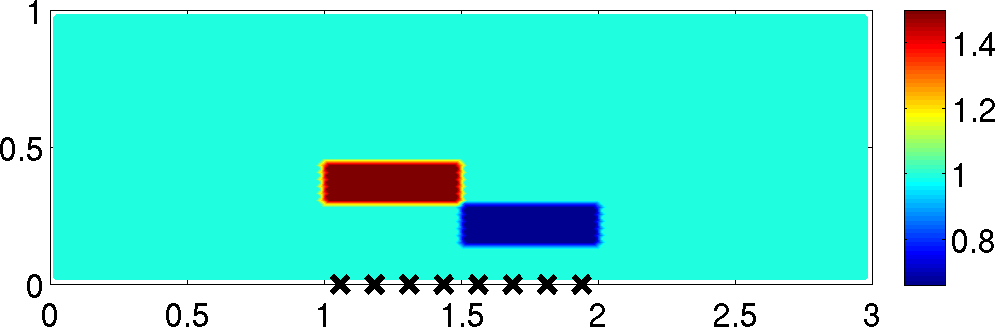}
\includegraphics[width=0.49\textwidth]{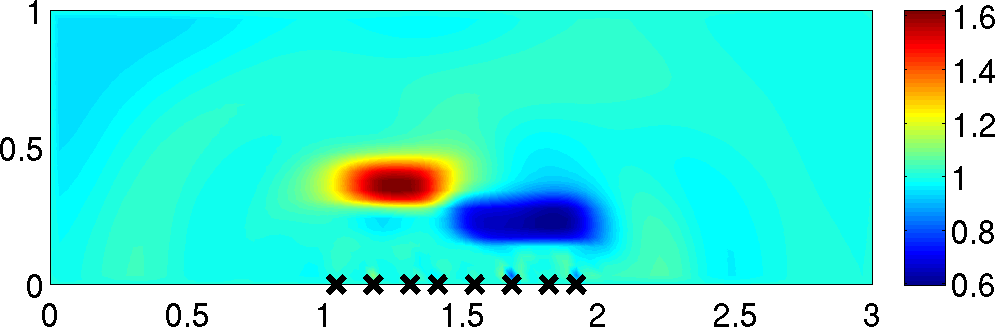}
\includegraphics[width=0.49\textwidth]{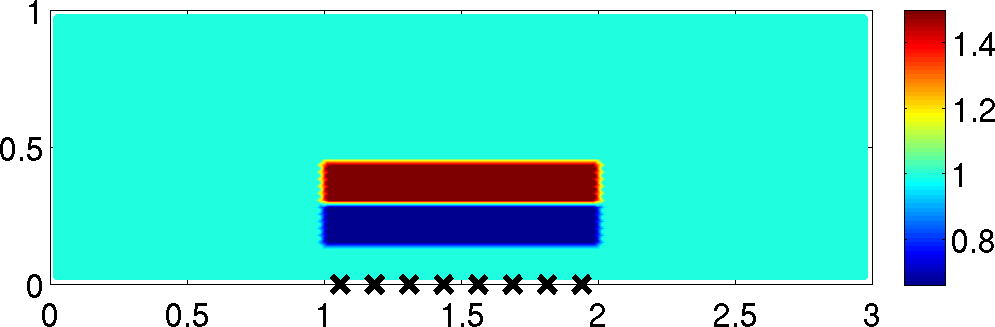}
\includegraphics[width=0.49\textwidth]{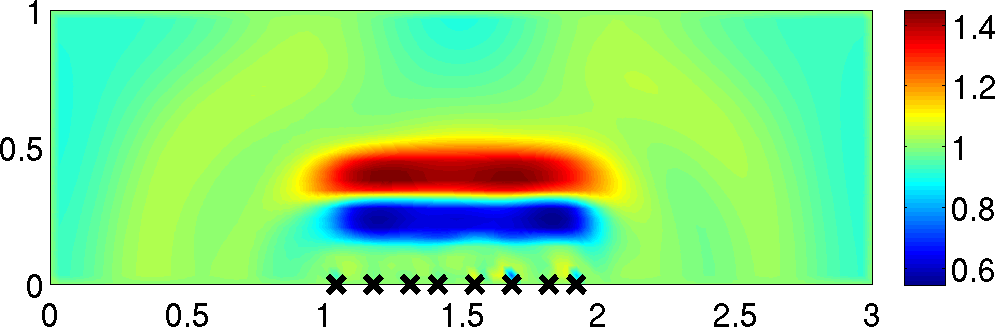}
\includegraphics[width=0.49\textwidth]{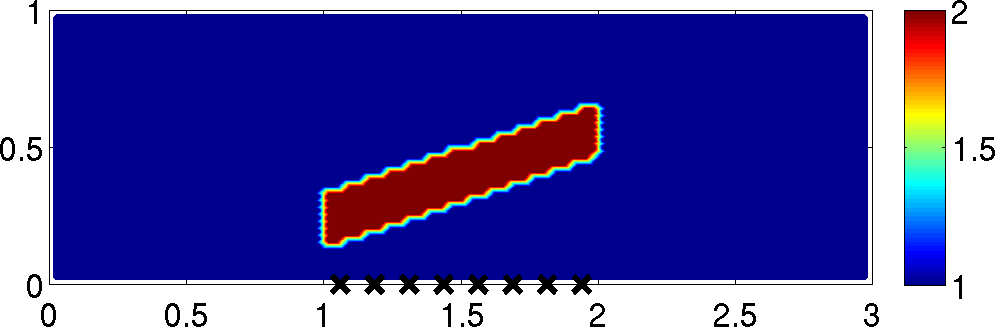}
\includegraphics[width=0.49\textwidth]{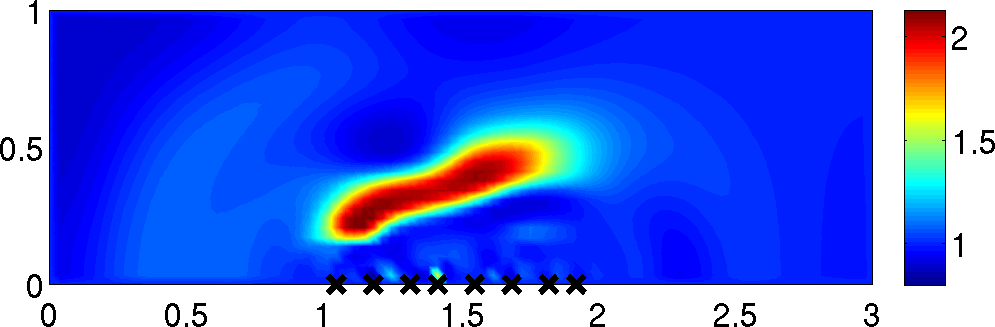}
\caption{Reconstructions in two dimensions in the rectangular domain
$\Omega = [0,3] \times [0,1]$. Left column: true coefficient $r(x)$; 
right column: reconstruction after a single Gauss-Newton iteration. 
Mid-points of sources/receivers $\cJ_j$ are marked with black $\times$,
$j=1,\ldots,N_d$, $N_d = 8$.
}
\label{fig:2Drec}
\end{center}
\end{figure}

We consider a two dimensional example in a rectangular domain 
$\Omega = [0,3] \times [0,1]$. The fine grid to simulate the data has the
dimension $120 \times 40$ nodes, while the coarse grid used in inversion
is $90 \times 30$ nodes. We use $N_d = 8$ sources/receivers with 
disjoint supports $\cJ_j$ uniformly distributed on the accessible boundary
interval $\cB_A = \{ (x_1,x_2) \;|\; x_1 \in (1,2), \; x_2 = 0 \}$. For each 
diagonal entry of the measured data matrix $y_{jj}(t)$, 
$j=1,\ldots,N_d$ a reduced order model with $m = 5$ is constructed.

As mentioned in section \ref{sec:2D}, the interpolation node 
$\widetilde{s}$ for the matching conditions (\ref{eqn:2Dmatchcond})
is chosen so that the sensitivity waves reach the boundary $\cB$
without reflecting from it. This is shown in Figure \ref{fig:2Dsens}
where the sensitivities for one particular source/receiver are 
plotted. For this $\Omega$ the interpolation node that gives the desired 
behavior is $\widehat{s} = 60$. We cannot take a smaller $\widehat{s}$
since the sensitivity function for $\kappa_5^{(4,4)}$ already touches
the boundary $x_2=1$. On the other hand, increasing $\widehat{s}$ will
shrink the region covered by the sensitivity functions and thus will
reduce the resolution away from $\cB_A$.

We solve the optimization problem (\ref{eqn:roptim2D}) for the two 
dimensional media using a regularized preconditioned Gauss-Newton 
inversion algorithm from section \ref{subsec:invalg} adapted for the 
objective functional of (\ref{eqn:roptim2D}). In two dimensions the 
preconditioner appears to be even more efficient with high quality 
reconstructions obtained after a single iteration. Subsequent iterations
improve the reconstruction marginally, so in Figure \ref{fig:2Drec}
we show the solutions after a single Gauss-Newton iteration starting
from a uniform initial guess $r(x) \equiv 1$.

All three examples in Figure \ref{fig:2Drec} are piecewise constant
so the inversion is regularized with a discrete $H^1$ seminorm ($W=I$).
In the first two examples there are two rectangular inclusions in each 
with the contrast from 
$\underset{x \in \Omega}{\mbox{min }} r(x) = 0.66$ to 
$\underset{x \in \Omega}{\mbox{max }} r(x) = 1.5$ on a unit 
background. The inclusions touch each other at a corner and a side
respectively. This demonstrates that the method handles well the sharp 
interfaces. In the third example there is a single tilted inclusion 
of contrast $2$ on a unit background. It is used to show the gradual 
loss of resolution away from $\cB_A$. All three examples are narrow 
aperture, meaning that the horizontal extent of the inclusions is equal 
to the width of $\cB_A$. Overall the reconstruction quality is good
with the contrast captured fully by the first Gauss-Newton iteration.
One can iterate further to improve the reconstruction, but an adaptive 
choice of the step length $\alpha^{(p)}$ is required for convergence.

\section{Summary}
\label{sec:conclusion}

We introduced a numerical inversion algorithm for linear parabolic
partial differential equations.  The problem arises in the application
of controlled source electromagnetic inversion, where the unknown is
the subsurface electrical resistivity $r(x)$ in the earth. We study
the inversion method in one and two dimensional media, but extensions
to three dimensions are possible.

To motivate the inversion algorithm we place the inverse problem in a
model reduction framework. We semidiscretize in $x$ the parabolic
partial differential equation on a grid with $N \gg 1$ points, and
obtain a dynamical system with transfer function $Y(s;\br)$, the
Laplace transform of the time measurements. In two dimensions the
transfer function is matrix valued, and we construct reduced models
separately, for each entry on its diagonal. Each model reduction
construction is as in the one dimensional case.

The reduced models are dynamical systems of much smaller size 
$m \ll N$, with transfer function $Y_m(s) \approx Y(s;\br)$. 
Because $Y_m(s)$ is a rational function of $s$, we solve a rational 
approximation problem.  We study various such approximants to 
determine which are best suited for inversion. We end up with a 
multipoint Pad\'{e} approximant $Y_m(s)$, which interpolates $Y$ 
and its first derivatives at nodes distributed geometrically in 
$\mathbb{R}^+$.  The inversion algorithm is a Gauss-Newton 
iteration for an optimization problem preconditioned with 
non-linear mappings $\cQ$ and $\cR$. These mappings are the 
essential ingredients in the inversion.

Most inversion algorithms estimate $r(x)$ by minimizing over
discretized resistivity vectors $\br \in \mathbb{R}_+^N$ the least
squares misfit of the data $d(t)$ and the mathematical model 
$y(t;\br)$ of the measurements.  By construction, our mapping $\cR(\br) =
\cQ(y(\cdot;\br))$ is an approximate identity when restricted to a
subset of sufficiently regular resistivities $\br$. We use it as a
non-linear preconditioner in the inversion, meaning that we minimize 
the least squares misfit of $\cQ(d(\,\cdot\,))$ and $\cR(\br)$. The 
advantage is the stability of the inversion and very fast convergence 
of the iteration.

We define the non-linear preconditioner $\cR(\br)$ via an explicit 
chain of non-linear mappings. Each step in the chain involves a 
numerically stable computation. The computation of the Jacobian 
$\cD \cR$ follows by the chain rule and we describe it explicitly, 
step by step. The only unstable computation in the inversion is 
the data fitting calculation of $\cQ(d(\,\cdot\,))$. The instability 
is inherited from that of the inverse problem and is unavoidable. We
mitigate it by restricting the size $m$ of the reduced model
adaptively, depending on the noise level. The smaller $m$ is, the
lower the resolution of the estimated resistivity. This is because at
each iteration the resistivity updates are in the range of $(\cD
\cR)^\dagger$, of low dimension $2 m \ll N$. We improve the results by
adding corrections in the null space of $\cD \cR$, so that we minimize
a regularization functional that incorporates prior information about
the unknown resistivity.  The performance of the algorithm is assessed
with numerical simulations in one and two dimensions.

\section*{Acknowledgments}
\label{sec:acknowledge}
The work of L. Borcea was partially supported by the AFSOR Grant
FA9550-12-1-0117, the ONR Grant N000141410077, and by the NSF Grant
DMS-0934594. The work of A. Mamonov was partially
supported by the NSF grants DMS-0934594, DMS-0914465 and DMS-0914840.

\appendix

\section{Computation of the non-linear preconditioner and its Jacobian}
\label{app:compjac}

The computation of the non-linear preconditioner 
$\cR(\br)= \cQ(y(\,\cdot\,;\br))$ and its Jacobian is the
most complex and time consuming computation in our inversion
scheme. Nevertheless, it is much more efficient than the traditional
inversion approach because it avoids the repeated computation of the
time domain solution of the forward problem.  We explain here the
details of the computation of $\cR$ and $\cD \cR$ via the chain of
mappings (\ref{eqn:kappachain}). We do so only for the multipoint
Pad\'{e} approximant, which we showed in section \ref{subsec:ratinterp} 
to be best suited for inversion.

\textbf{(a)~} The matrix $A$ is defined by (\ref{eqn:Ar}) for a given
$\br$.  Differentiating (\ref{eqn:Ar}) yields
\begin{equation}
\frac{\partial A}{\partial r_k} = - D^T \be_k \be_k^T D = - \bd_k \bd_k^T,
\quad k=1,\ldots,N,
\end{equation} 
with $\bd_k = D_{k,\,1:N}^T$. This is a rank one matrix.

\textbf{(b)~} At this step we differentiate the orthonormal basis 
$V$ of the Krylov subspace $\cK_m (\widetilde{\bf s})$. There are 
different ways of computing an orthonormal basis of 
$\cK_m (\widetilde{\bf s})$. One choice is to use a rational 
Lanczos algorithm \cite{ratlanczos}. While this may be a more stable way 
of computing $V$ compared to other approaches, differentiation formulas 
are difficult to derive and implement.  We consider an alternative approach, 
based on the differentiation of the QR decomposition.

We compute first the matrix 
\begin{equation*}
K = \left[ (\widetilde{s}_1 I - A)^{-1} \bb, \ldots,
(\widetilde{s}_m I - A)^{-1} \bb \right] \in 
\mathbb{R}^{N \times m},
\end{equation*}
and its derivatives
\begin{equation*}
\frac{\partial K}{\partial r_k} = -\left[ (\widetilde{s}_1 I -
  A)^{-1} \bd_k [(\widetilde{s}_1 I - A)^{-1} \bd_k]^T \bb, \ldots,
  (\widetilde{s}_m I - A)^{-1} \bd_k [(\widetilde{s}_m I -
    A)^{-1} \bd_k]^T \bb \right].
\end{equation*}
Then $V$ can be defined via the QR decomposition of $K$, which we
write as
\begin{equation*}
K = V U,
\label{eqn:KQR}
\end{equation*}
where $V \in \mathbb{R}^{N \times m}$ is orthogonal $V^T V = I_m \in \mathbb{R}^{m
  \times m}$, and $U \in \mathbb{R}^{m \times m}$ is upper
triangular. If we denote $L = U^T$, then (\ref{eqn:KQR}) implies
\begin{equation*}
K^T K = L L^T.
\end{equation*}
That is to say $L = U^T$ is a Cholesky factor of $K^T K$. At the same
time, when we differentiate (\ref{eqn:KQR}), we obtain
\begin{equation}
\frac{\partial V}{\partial r_k} = \left( \frac{\partial K}{\partial r_k} - 
V \frac{\partial U}{\partial r_k} \right) U^{-1}, \quad k=1,\ldots,N.
\label{eqn:vdiff}
\end{equation}
Since we already know $\partial K / \partial r_k$, it remains to
compute the derivative $\partial U / \partial r_k$ of the Cholesky
factorization of $K^T K$. This is given in the following proposition,
which can be proved by direct computation once we write
$
\delta (L L^T) = (\delta L) L^T + L (\delta L)^T
$
and solve for the columns of $\delta L$, one a time, using that
$\delta L$ is lower triangular.
\begin{prop}[Differentiation of Cholesky factorization] \label{prop:choldiff}
Let $M \in \mathbb{R}^{m \times m}$ be a matrix with Cholesky
factorization $M = L L^T$.  Given the perturbation $\delta M$ of $M$,
the corresponding perturbation $\delta L$ of the Cholesky factor is
computed by the following algorithm.
\begin{itemize}
\item[] For $k=1,\ldots,m$
\item[] \hskip0.1\textwidth $\delta L_{kk} = \dfrac{1}{L_{kk}} \left( 
\dfrac{\delta M_{kk}}{2} - \sum\limits_{j=1}^{k-1} \delta L_{kj} L_{kj} \right).$
\item[] \hskip0.1\textwidth For $i=k+1,\ldots,m$
\item[] \hskip0.2\textwidth $\delta L_{ik} = \dfrac{1}{L_{kk}} \left( \delta M_{ik}
- \sum\limits_{j=1}^{k} \delta L_{kj} L_{ij}
- \sum\limits_{j=1}^{k-1} \delta L_{ij} L_{kj}
\right). $
\end{itemize}
\end{prop}

\noindent We use Proposition \ref{prop:choldiff} for $M = K^T K = L L^T$, with
perturbation
\[ \delta M = \frac{\partial M}{\partial r_k} \delta r_k, \qquad
\frac{\partial M }{\partial r_k} = \left(\frac{\partial K }{ \partial
  r_k}\right)^T K + K^T \left(\frac{\partial K }{ \partial
  r_k}\right),\] to obtain 
\[
\delta L = \frac{\partial L}{\partial r_k} \delta r_k \quad \mbox{and} 
\quad 
\frac{\partial U }{ \partial r_k} = \left( \frac{\partial L}{\partial
  r_k}\right)^T.\] The computation of $\partial V / \partial r_k$
follows from (\ref{eqn:vdiff}).

\textbf{(c)~} Once we have $V$ and its derivatives we compute from 
(\ref{eqn:rom-project})
\begin{eqnarray}
\frac{\partial A_m}{\partial r_k} & = & - V^T \bd_k \bd_k^T V +
\frac{\partial V}{\partial r_k}^T A V + V^T A \frac{\partial
  V}{\partial r_k}, \\ \frac{\partial \bb_m}{\partial r_k} & = &
\frac{\partial V}{\partial r_k}^T \bb.
\end{eqnarray}

\textbf{(d)--(e)~} There are two possible ways to go from the reduced
model $A_m$, $\bb_m$ to the continued fraction coefficients $\kappa_j$,
$\widehat{\kappa}_j$. Both approaches use a Lanczos iteration to obtain
a symmetric tridiagonal matrix, which we denote by
\begin{equation}
T = \begin{bmatrix}
\alpha_1 & \beta_2  &         & \\
\beta_2  & \alpha_2 & \ddots  & \\
         & \ddots   & \ddots  & \beta_m \\
         &          & \beta_m & \alpha_m
\end{bmatrix}.
\end{equation}

Note that for an arbitray symmetric matrix $E \in \mathbb{R}^{m \times
  m}$ we can compute a tridiagonal matrix $T$ that is unitarily
similar to $E$ via a Lanczos iteration with full
reorthogonalization. Since the dimension $m$ of the problem is small,
we reorthogonalize at every step to ensure maximum numerical
stability.  The iteration is as follows:

\vspace{0.1in}
\begin{itemize}
\item[] Initialize $\bx_1 = \bfeta$, $\bx_0=0$, $\beta_1 = 0$.
\item[] For $j=1,\ldots,m-1$
\item[] \hskip0.1\textwidth $\alpha_j = \bx_j^T E \bx_j$,
\item[] \hskip0.1\textwidth $\widetilde{\bu}_{j+1} = E \bx_j - \alpha_j
  \bx_j - \beta_j \bx_{j-1}$,
\item[] \hskip0.1\textwidth $\bu_{j+1} = (I - X_{1:N,\,1:j}
  (X_{1:N,\,1:j})^T) \widetilde{\bu}_{j+1}$,
\item[] \hskip0.1\textwidth $\beta_{j+1} = \| \bu_{j+1} \|$,
\item[] \hskip0.1\textwidth $\bx_{j+1} = \dfrac{\bu_{j+1}}{\beta_{j+1}}$.
\item[] $\alpha_m = \bx_m^T E \bx_m$,
\end{itemize}

\vspace{0.1in} \noindent We have two choices of the initial vector
$\bfeta$ and the matrix $E$ such that $T = X^T E X$ with $ X^T X =
I_m$. The first takes $E = A_m$, and $\bfeta = \bb_m / \| \bb_m \|$.
It combines steps (e) and (d) and goes from the reduced model $A_m$,
$\bb_m$ directly to the tridiagonal matrix $T$.  The second approach
is to compute first the eigenvalue decomposition (\ref{eqn:thetac}) of
$A_m$ to get the poles $-\theta_j$ and the residues $c_j$, for
$j=1,\ldots,m$. Then take $E = - \mbox{diag}(\theta_1, \ldots,
\theta_m)$ and 
\begin{equation}
\bfeta = (\eta_1, \ldots, \eta_m)^T, \qquad 
\eta_i = \sqrt{\frac{c_i}{\sum\limits_{s=1}^{m} c_s}}, \quad i=1,\ldots,m.
\label{eqn:etac}
\end{equation}

From the computed $T$ with either of the two approaches, we obtain the
coefficients of the continued fraction using the formulas from
\cite{druskin2000gaussian}.
\begin{eqnarray}
\widehat{\kappa}_1 &=& \frac{1}{\sum\limits_{s=1}^{m} c_s}, \quad
\kappa_1 = - \frac{1}{ \widehat{\kappa}_1
  \alpha_1}, \label{eqn:kappa-alpha-beta-1} \\ \widehat{\kappa}_j & =
& \frac{1}{ \kappa_{j-1}^2 \beta_{j-1}^2 \widehat{\kappa}_{j-1} },
\quad j=2,\ldots,m, \\ \kappa_j & = & - \frac{1}{ \alpha_j
  \widehat{\kappa}_j + \dfrac{1}{\kappa_{j-1}} }, \quad j=2,\ldots,m.
\label{eqn:kappahat-alpha-beta-j}
\end{eqnarray}

To differentiate the mappings (d) and (e), we need to differentiate
the Lanczos iteration. In general, there is no explicit formula for
the perturbations of the entries of $T$ in terms of perturbations of
$E$ and $\bfeta$.  Thus, we would need to differentiate the Lanczos
iteration directly, with an algorithm that computes the perturbations
of $\alpha_j$, $\beta_j$ and $\bx_j$ iteratively, for increasing
$j$. However, for the second approach described above, we can use the
explicit perturbation formulas for the Lanczos iteration derived in
\cite{borcea2005continuum}. To apply these formulas, we need to
differentiate the steps (d) and (e) separately.

The differentiation of (d) follows directly from the differentiation
of the eigendecomposition (\ref{eqn:thetac}) of $A_m$.  It is given by
\begin{eqnarray}
\frac{\partial \theta_j}{\partial r_k} & = & - \bz_j \frac{\partial
  A_m}{\partial r_k} \bz_j, \\ \frac{\partial \bz_j}{\partial r_k} & = & -
(A_m + \theta_j I)^\dagger \frac{\partial A_m}{\partial r_k} \bz_j,
\\ \frac{\partial c_j}{\partial r_k} & = & 2 (\bb_m^T \bz_j) \left( \bb^T
\frac{\partial V}{\partial r_k} \bz_j + \bb_m^T \frac{\partial
  \bz_j}{\partial r_k} \right),
\label{eqn:diffc}
\end{eqnarray}
where $\dagger$ denotes the pseudoinverse.

For the computation of the derivatives in step (e) we use the explicit
formulas from \cite{borcea2005continuum}, for $E = \Theta$ and $X =
Q^T$ and $\bfeta$ given by (\ref{eqn:etac}).  Let us define the
vectors
\begin{equation}
\bfdelta_\alpha = \begin{pmatrix} \delta \alpha_1 \\ \delta \alpha_1 +
  \delta \alpha_2 \\ \vdots \\ \sum\limits_{j=1}^{m} \delta\alpha_j
\end{pmatrix}, \quad
\bfdelta_\beta = \begin{pmatrix} \delta\beta_1 / \beta_1 \\ \delta
  \beta_1 / \beta_1 + \delta \beta_2 / \beta_2 \\ \vdots
  \\ \sum\limits_{j=1}^{m} \delta\beta_j / \beta_j
\end{pmatrix}.
\end{equation}
They can be expressed in terms of the perturbations of the poles
$\bftheta = (\theta_1, \ldots,\theta_m)^T$ and the initial vector $\bfeta$ as
\begin{equation*}
\bfdelta_\alpha =  - A_\theta \delta \bftheta + A_\eta \delta
\bfeta, \qquad \bfdelta_\beta =  - B_\theta \delta \bftheta + B_\eta
\delta \bfeta.
\end{equation*}
Here $A_\theta, A_\eta \in \mathbb{R}^{m
  \times m}$ are the matrices with entries 
\begin{eqnarray*}
A_\theta^{ij} & = & 1 + \beta_i
\mathop{\sum\limits_{p=1}^{m}}\limits_{p \neq j} \frac{1}{\theta_p -
  \theta_j} \left[ 2 Q_{ip} Q_{i+1,p} - \frac{Q_{1 p}}{Q_{1 j}}(Q_{i
    p} Q_{i+1,j} - Q_{i+1,p} Q_{i j}) \right], \\ 
A_\eta^{ij} & = & 2 \beta_i \frac{Q_{i+1, j} Q_{i j}}{Q_{1 j}},
\quad A_\eta^{m j} = 0, \quad A_\theta^{mj} = 1,
\end{eqnarray*}
for $i=1,\ldots,m-1$, $j=1,\ldots,m$ and the entries of $B_\theta,
B_\eta \in \mathbb{R}^{(m-1) \times m}$ are
\begin{eqnarray*}
B_\theta^{ij} & = & \mathop{\sum\limits_{p=1}^{m}}\limits_{p \neq j}
\frac{1}{\theta_p - \theta_j} \left[ (Q_{i+1,p})^2 - \frac{Q_{1
      p}}{Q_{1 j}} Q_{i+1,p} Q_{i+1,j} \right], \\ B_\eta^{ij} & = &
\frac{(Q_{i+1, j})^2}{Q_{1 j}},
\end{eqnarray*}
for $i=1,\ldots,m-1$ and $j=1,\ldots,m$.

Note that the computation of $\delta \bfeta$ can be done by
differentiating (\ref{eqn:etac}) and using (\ref{eqn:diffc}). Once the
vectors $\bfdelta_\alpha$ and $\bfdelta_\beta$ are known, it is trivial to
obtain the individual perturbations $\delta \alpha_j$ and $\delta
\beta_j$. From those perturbations we compute the derivatives
of $\kappa_j$ and $\widehat{\kappa}_j$ by differentiating relations
(\ref{eqn:kappa-alpha-beta-1})--(\ref{eqn:kappahat-alpha-beta-j}).
The computation is straightforward and we do not include it here.

\textbf{(f)} The differentiation at the last step in the chain of
mappings (\ref{eqn:kappachain}) is trivial. It is just the
derivative of the logarithm.

\bibliography{biblio} 
\bibliographystyle{siam}
\end{document}